\documentclass[12pt]{amsart}   
\usepackage{amsmath}
\usepackage{amssymb}

\newtheorem{theorem}{Theorem}[section]   
\newtheorem{lemma}[theorem]{Lemma}   
\newtheorem{corollary}[theorem]{Corollary}   
\newtheorem{proposition}[theorem]{Proposition}   
   
\theoremstyle{definition}

\newtheorem{definition}{Definition}

\numberwithin{equation}{section}   
   
\renewcommand{\Re}{\mathrm{Re}\,}
\renewcommand{\Im}{\mathrm{Im}\,}

\def\prt{\partial} 
\newcommand{\Lep}{{L_\epsilon^+}}
   
\def\reals{ {{\mathbb R}} }   
   
\def\integers{ {{\mathbb Z}} }   
\def\complex{{\mathbb C}}

\font\cal=cmsy10 scaled\magstep1

\def\hils{\text{\cal H}}
\def\Z{\integers}
\def\R{\reals}
\def\C{\complex}
  
\newcommand\beq{\begin{equation}}   
\newcommand\eeq{\end{equation}}   

\def\tr{\mathrm{Tr}\,}  

\newcommand{\floor}[1]{\lfloor#1\rfloor}

\begin{document}


\author{Rowan Killip}   
\address{Department of Mathematics \\
University of Pennsylvania \\   
Philadelphia, PA 19104-6395, USA}   
\email{killip@math.upenn.edu} 
\thanks{RK was supported in part by NSF grant No.~DMS-9729992}   
\author{Alexander Kiselev}   
\address{Department of Mathematics \\   
University of Chicago \\   
Chicago, IL 60637, USA}   
\email{kiselev@math.uchicago.edu}   
\thanks{AK was supported in part by NSF grant No.~DMS-0102554, and by an Alfred P. Sloan Research Fellowship}   
\author{Yoram Last}
\address{Institute of Mathematics, The Hebrew University \\
91904 Jerusalem, Israel}
\email{ylast@math.huji.ac.il}
\thanks{YL was supported in part by NSF grant No.~DMS-9801474, by The Israel Science Foundation 
Grant No.~447/99, and by an Allon fellowship.}   

\date{December 8, 2001}   
   
\title[bounds on wavepacket spreading]   
{Dynamical upper bounds on wavepacket spreading}   
   
\begin{abstract}  
We derive a general upper bound on the spreading
rate of wavepackets in the framework of Schr\"odinger
time evolution. Our result consists of showing that
a portion of the wavepacket cannot escape outside a
ball whose size grows dynamically in time, where
the rate of this growth is determined by
properties of the spectral measure and by
spatial properties of solutions
of an associated time independent Schr\"odinger
equation. We also derive a new lower bound on
the spreading rate, which is strongly connected
with our upper bound. We apply these new bounds to
the Fibonacci Hamiltonian---the most studied
one-dimensional model of quasicrystals. As a result, we
obtain for this model upper and lower dynamical bounds
establishing wavepacket spreading rates which are
intermediate between ballistic transport and localization. 
The bounds have the same qualitative behavior in the
limit of large coupling. 
\end{abstract}   

\maketitle   
   
\section{Introduction} 

Let $H$ be a self adjoint operator on a separable Hilbert space $\hils$.
The time-dependent Schr\"odinger equation of quantum mechanics,
$i\partial_t\psi = H\psi$, leads to a unitary
dynamical evolution in $\hils$, given by
\begin{equation}\label{evolution}
\psi(t)=e^{-iHt}\psi.
\end{equation}
We are primarily interested here in cases where $H$ is a Schr\"odinger
operator of the form $\Delta + V$ on a space such as $\ell^2(\Z^d)$ or
$L^2(\R^d)$ and where the initial vector $\psi$ is
a localized wavepacket. Under the evolution \eqref{evolution},
the wavepacket will often spread with time and it is the nature
of this spreading which is our main object of interest here. This question,
``What determines the spreading of a wavepacket'' \cite{KKKG}, has been
an active field of research over the last 15 years and there is by now
a considerable body of literature devoted to it (\cite{BCM,BGT,BSb,BT,Combes,CM,
DJLS,GdB,Guarneri,Gu2,GM,GSb1,GSb2,KKKG,KL,Last,Mantica1,Mantica2,Piechon,Simon1,
Simon2,Tch1,Tch2,WA} is just a partial list).
We note that there are some situations, such as some systems with
absolutely continuous spectrum which can be studied by scattering theory
\cite{RS3}, where good understanding of the dynamics exists and is much
older. The primary focus of the more recent studies has been on greater
generality and, in particular, on systems exhibiting spectral
phenomena such as singular continuous spectrum and dense point spectrum.
In particular, singular continuous spectrum tends to occur for basic
models of quasicrystals (see, e.g., \cite{Dam1,Suto3}) and the study
of wavepacket spreading in such models is important to the understanding
of transport properties in such media.

A major focus of many recent studies has been to elucidate the relations
between wavepacket spreading rates and continuity properties of spectral
measures. Recall that each initial vector $\psi$ has a spectral measure
$\mu_\psi$, defined as the unique Borel measure obeying
$
\langle f(H)\psi,\psi\rangle = \int_{\sigma (H)} f(x)\,d\mu_\psi(x)
$
for any measurable function $f$. Here we use $\langle\cdot, \cdot\rangle$
to denote the scalar product in $\hils$. A major discovery of Guarneri
\cite{Guarneri,Gu2}, which has also been extended by several other authors
\cite{BGT,BT,Combes,GSb1,GSb2,Last},
is that appropriately measured continuity properties of the spectral measure
$\mu_\psi$ imply fairly general lower bounds on the spreading rate 
of the wavepacket. Recently, Kiselev-Last \cite{KL} (also see \cite{KKKG}
for a related heuristic result and \cite{Tch2} for a recent extension)
proved a variant of Guarneri's bound which blends continuity properties of
the spectral measure with the spatial decay rate of generalized eigenfunctions.
This bound is generally stronger than what one can obtain from properties
of the spectral measure alone and it is also of somewhat broader applicability.
We note that while there is some pure theoretical interest in relating
wavepacket dynamics to such things as continuity properties of
spectral measures and spatial properties of eigenfunctions, the interest
in the above bounds is more than purely theoretical. Properties of
spectral measures and eigenfunctions can be studied for concrete
models of interest and the above discussed bounds have been used to obtain
dynamical results in cases where there is currently no alternative
approach to study the dynamics. In particular, this approach has been
used to obtain dynamical results for a number of quasiperiodic operators 
\cite{Dam,DKL,JL2} and for operators with decaying potentials \cite{KL}.

The situation with putting upper bounds on wavepacket spreading rates
is much more problematic. There is a fairly general argument of Simon
\cite{Simon1} showing that for a broad class of Schr\"odinger operators
(in particular, every operator of the form $\Delta+V$ on $\ell^2(\Z^d)$)
wavepacket spreading rates cannot be more than ballistic (namely,
linear in time). For the regime of dynamical localization (namely,
situations where the wavepacket does not spread beyond a finite region
of space and so, in particular, the spectrum must be pure point) there
is a fairly good understanding \cite{DJLS,GdB,Tch1} of how suitable spatial
properties of eigenvectors (or of spectral projections) are related to dynamics and
how to specify necessary and sufficient conditions for the occurrence
of complete dynamical localization. For problems with continuous spectrum,
however, it had been a noted open problem to get any results beyond
Simon's ballistic bound, namely, to be able to specify conditions
that would guarantee wavepacket spreading rates that are slower than
ballistic.

The main aim of the current paper is to provide a certain solution to
this problem. 
What we are
able to do, in essence, is to specify conditions that would guarantee
that {\it some} portion of the wavepacket does not spread too
fast (namely, it remains within a box whose size grows with time
at an appropriate rate). We cannot control the entire wavepacket
and thus we cannot bound the growth rate of quantities that are
sensitive to small, fast spreading portions of the wavepacket,
such as moments of the position operator. Nevertheless, we believe
that our result is in many ways a natural complement to some of the Guarneri
type lower bounds discussed above. We note that while Guarneri type
bounds, such as the variants obtained in \cite{KL,Last}, have been
often formulated in terms of moments of the position operator, such bounds
are usually obtained by controlling only a portion of the wavepacket.
It so happens that in order to establish a {\it lower} bound on the
growth rate of moments of the position operator it suffices to
show that {\it some} portion of the wavepacket is spreading at an
appropriate rate, whereas bounding such quantities from above would require
tight control of the entire wavepacket---a much more delicate task.
Moreover, we will see that our upper bound has a strong connection with
the Kiselev-Last \cite{KL} lower bound discussed above. While we formulate
things somewhat differently here, the basic information going into
our upper bound is, roughly speaking, the same combination of spectral
measure continuity and spatial behavior of generalized eigenfunctions which goes
into that lower bound. Indeed, we also formulate here a new stronger variant
of this lower bound showing that, in essence, we have good understanding
of how the spreading rate of the slowest moving portion of the wavepacket
is determined. The behavior of the fastest moving portion remains outside
our scope of control and, in particular, establishing sub-ballistic upper
bounds on the growth rate of moments of the position operator remains
an interesting open problem.

{From} here on we study tridiagonal operators of the form 
\begin{equation}\label{jacobi}
(H u)(n)= a(n) u(n+1) + a (n-1) u(n-1) + b (n) u(n)
\end{equation}
(also called Jacobi matrices) defined on $\ell^2(\integers_+)$ or 
$\ell^2(\integers)$. Here $b(n)$ and $a(n)$ are real numbers, and 
$a(n) \ne 0$ for any $n$. Moreover, we assume that $\sum_{n =\pm 1}^{\pm \infty}
|a(n)|^{-1} = \infty,$ which is sufficient to ensure that these operators are
self-adjoint (see, e.g., \cite{Ber}). We note that 
discrete one-dimensional Schr\"odinger operators of the form $\Delta+V$ on $\ell^2(\integers_+)$ or 
$\ell^2(\integers)$ are a special case of such tridiagonal operators
(obtained by setting $a(n) \equiv 1$ and $b(n)=V(n)$). While such one-dimensional
operators will be our main interest here, we discuss in Section 7 below how
our results are also applicable to more general problems and, in particular, how
our dynamical upper bound is also applicable to multidimensional Schr\"odinger
operators of the form $\Delta+V$ on $\ell^2(\Z^d)$.

We will look at the evolution of the initial vector 
$\psi = \delta_1$ (where $\delta_n(k) = \delta_{nk}$),
but the results can also be recast for other $\ell^2$ 
vectors.

Let $u_\theta(n,z)$ be the solution of the equation 
\begin{equation}\label{ef}
a(n)u(n+1) + a(n-1) u(n-1) + b(n) u(n) = z u(n),
\end{equation}
$z \in \complex,$
satisfying $u_\theta(0,z)=\sin \theta,$ $u_\theta(1,z) = \cos \theta.$ 
Given a function $\phi:\Z_+\to\C,$ 
we define for any $L\geq 0$,
\begin{equation}\label{normL}
 \|\phi\|^2_L = \sum\limits_{n=1}^{\floor{L}} |\phi(n)|^2 + (L-\floor{L})|\phi(\floor{L}+1)|^2 ,
\end{equation}
where we use the convention that the sum is zero if the upper limit
is less than the lower one.
Similarly, for a function $\phi:\Z\to\C$, we define
\[
  \|\phi\|^2_{L_1,L_2} = \sum\limits_{n=-\floor{L_1}}^{\floor{L_2}} |\phi(n)|^2 +
(L_1-\floor{L_1})|\phi(-\floor{L_1}-1)|^2+(L_2-\floor{L_2})|\phi(\floor{L_2}+1)|^2 .
\]
Here $\floor{x}$ 
is the greatest integer less than or equal to $x$.
While the extension to real $L$ is not essential, it will be technically convenient.  
Consider first the half-line problem, namely, the operator on $\ell^2(\Z_+)$.
Without loss of generality, we restrict consideration to a Dirichlet boundary
condition ($\theta =0$), where the operator \eqref{jacobi}
is defined on functions vanishing for $n \leq 0.$ We denote the resulting self-adjoint operator 
by $H_+.$ We also assume here $a(0)=1$. Note that the operator $H_+$ itself is independent
of $a(0)$. The choice of $a(0)$ is only effecting the definition of the solutions
$u_\theta(n,z)$ for $\theta\not= 0$.
Given $\epsilon >0$ and an energy $E \in \reals,$ we define the length scale $L^+_\epsilon(E)$ by 
\begin{equation}\label{ls}
\max_\theta \|u_\theta (n,E)\|_{L^+_\epsilon (E)}\,\cdot\, \min_\theta \|u_\theta (n,E)\|_{L^+_\epsilon (E)}
= \epsilon^{-1}.
\end{equation}
The left-hand side is monotonely increasing as the region 
of summation expands and so $L^+_\epsilon(E)$ is well-defined for every 
$\epsilon > 0$. We have $L^+_\epsilon(E)\to 1$ as $\epsilon\to\infty$
and $L^+_\epsilon(E)\to\infty$ as $\epsilon\to 0$.
Let $m_+(z)$ be the Weyl-Titchmarsh $m$ function corresponding to the operator $H_+,$
\[ m_+(z) = \langle (H_+ -z)^{-1}\delta_1, \delta_1 \rangle \] 
The spectral measure $\mu^+$ of the operator $H_+$ is defined by the equality
\[ m_+(z) = \int_\reals \frac{d \mu^+(x)}{x-z}, \]
and coincides with $\mu_{\delta_1}$, the spectral measure of the vector $\delta_1$.
As in previous works on the subject, we deal with time-averaged quantities.
For any function $A(t)$ of time, we denote 
\[ \langle A(t) \rangle_T = \tfrac{2}{T} \int\limits_0^\infty
e^{-2t/T} A(t) \, dt. \]
This type of averaging is slightly different from the Ces\`aro averaging
used in many previous works and it is more convenient for what we do here.
Note that the difference is not very significant.
In particular, power law behaviors in $T$ must be the same for both kinds of
average.

Our first result is the following 

\begin{theorem}\label{hld}
Let $H_+$ be the half-line operator defined on $\ell^2(\Z_+)$ by \eqref{jacobi}
and a Dirichlet boundary condition.
Let the characteristic scale $L^+_{T^{-1}}(E)$ be defined by \eqref{ls}.
Then for any $T>0$ and $L>1$, we have 
\begin{equation}\label{dynfir}
\langle \|e^{-iH_+ t} \delta_1 \|^2_{L} \rangle_T > C \mu^+\left(\left\{ E \left| \right. L^+_{T^{-1}}(E) \leq L \right\}\right),
\end{equation} 
where $C$ is a universal positive constant.
\end{theorem}
\noindent{\it Remarks.} 1. The expression on the left hand side of \eqref{dynfir}
gives the averaged norm of the portion of the wavepacket remaining
in a ball of size $L.$ Theorem~\ref{hld} shows that this norm
is bounded from below by a quantity proportional to the norm of the spectral projection of $\delta_1$ 
on the set of energies where $L^+_{T^{-1}}(E) \leq L$. Thus, if we choose
the size $L$ to be greater than a certain scale which depends on time 
and properties of the solutions $u_\theta$, we are guaranteed to have
a significant portion of the wavepacket remaining (on average) in a ball
of this size. \\
2. While the universal constant $C$ can be explicitly estimated from our proof
below, we made no real effort to obtain an optimal value for it. Our
technique can only yield a number which is significantly smaller than $1$ and
thus \eqref{dynfir} can control the location of only a portion of the
wavepacket. \\

The next criterion is a simple corollary of Theorem~\ref{hld}. It relates an upper bound on 
wavepacket spreading to the growth of the norms of transfer matrices.
Let us denote by $\Phi(n,E)$ the transfer matrix from the site $0$ 
to $n:$ 
\[
\Phi(n,E) = \left( \begin{array}{cc} u_0(n+1,E) & u_{\pi/2}(n+1,E) \\
u_{0}(n,E) & u_{\pi/2}(n,E) 
\end{array} \right).
\]
We also denote, for $L\geq 1$,
\[ 
\| \Phi(E) \|_L^2 = \sum\limits_{n=1}^{\floor{L}-1} \|\Phi(n,E)\|^2 + 
(L-\floor{L})\|\Phi(\floor{L},E)\|^2,
\]
where $\| \Phi(n,E) \|$ is the usual operator norm of the matrix $\Phi(n,E).$ 
Define $\tilde{L}_{\epsilon}^+ (E)$ by 
\begin{equation}\label{chartran}
 \|\Phi(E)\|_{\tilde{L}_{\epsilon}^+ (E)} = 2 \|\Phi(1,E)^{-1}\| \epsilon^{-1}.
\end{equation}
Then we have the same result for a new characteristic scale:

\begin{corollary}\label{hldt}
Let $H_+$ be the half-line operator defined on $\ell^2(\Z_+)$ by \eqref{jacobi}
and a Dirichlet boundary condition.
Let the characteristic scale $\tilde{L}^+_{T^{-1}}(E)$ be defined by \eqref{chartran}.
Then for any $T>0$ and $L\geq 2$, we have 
\begin{equation}
\langle \|e^{-iH_+ t} \delta_1 \|^2_{L} \rangle_T > C \mu^+\left( \left\{ E \left| \right. \tilde{L}^+_{T^{-1}}(E) \leq L \right\}\right),
\end{equation} 
where $C$ is a universal positive constant.
\end{corollary}
\noindent{\it Remark.}
By the Landauer formula \cite{SK}, the resistance
$\rho(n)$ of a sample of size $n$ is given by 
\[ \rho(n) = \tfrac{1}{2} \left[ \tfrac{1}{2} \tr (\Phi^t (n,E) \Phi(n,E))-1\right] 
\]
(where $\Phi^t$ is the transpose of $\Phi$).
Hence, the growth rate (with $L$) of $\| \Phi(E) \|_L^2$ is connected with
the growth rate of the sum
\[
\sum\limits_{n=1}^{L} \rho(n) .
\]
Corollary~\ref{hldt} can thus be viewed as a confirmation of a physically
reasonable statement that higher resistance leads to slower transport. \\

We will show below that $L^+_\epsilon(E) \leq \tilde{L}_{\epsilon}^+ (E),$ and hence
Corollary~\ref{hldt} is strictly weaker than Theorem~\ref{hld}. 
However, it may be more straightforward to apply it in some situations.

In addition, we derive a new lower bound on wavepacket spreading which is also
related to the characteristic length scale $L_\epsilon^+(E).$ 
\begin{theorem}\label{ldb}
Let $H_+$ be the half-line operator defined on $\ell^2(\integers_+)$
by \eqref{jacobi} with a Dirichlet boundary condition.
Denote by $P_S$ the spectral projection on some Borel set $S.$ 
Then for every $T>0$ and $L>0,$  
\begin{equation}\label{eqldb}
\langle \| e^{-iH_+t} P_S \delta_1 \|^2_{L}\rangle_T \leq C \int\limits_S \frac{\|u_0(n,E)\|^2_L}
{\|u_0(n,E)\|^2_{L_{T^{-1}}^+(E)}} d \mu^+(E),
\end{equation}
where $C$ is a universal constant. 
\end{theorem}
\noindent Roughly, Theorem~\ref{ldb} shows that if $L$ is such that the ratio  
\[ \|u_0(n,E)\|^2_L/
\|u_0(n,E)\|^2_{L_{T^{-1}}^+} \]
is small for $E \in S,$ then most of the portion of the wavepacket
which is supported on energies within the set $S$ leaves the
ball of radius $L$ by the time $T.$  

While we defined the scales $L_\epsilon^+(E)$ purely in terms of
solutions of the equation \eqref{ef}, we will see below that the quantities
which enter in this definition also determine the local scaling properties
of the spectral measure. This connection can be turned around to control
$L_\epsilon^+(E)$ using partial information on solution behavior along with
continuity properties of the spectral measure. In particular, 
the meaning of Theorem~\ref{ldb} and its connection with the Kiselev-Last
lower bound \cite{KL} can be clarified by noting that it implies the
following. 
Let $\alpha(E)$ be a measurable function such that for a.e. $E$ w.r.t. $\mu^+,$ 
\begin{equation}\label{ud}
 D^{\alpha(E)} \mu^+(E) = \limsup_{\epsilon \rightarrow 0}
 \frac{\mu^+(E+\epsilon, E-\epsilon)}{(2\epsilon)^{\alpha(E)}} < \infty.
\end{equation}
Let $\gamma(E)$ be a 
measurable function such that for a.e. $E$ w.r.t. $\mu^+,$ 
\begin{equation}\label{efb}
 \limsup_{L \rightarrow \infty} L^{-\gamma(E)} \| u_0 (n,E)\|^2_L < \infty. 
\end{equation}
We will show that Theorem~\ref{ldb} implies:
\begin{theorem}\label{pb}
Let $\alpha(E),$ $\gamma(E)$ satisfy \eqref{ud}, \eqref{efb} and let $\eta (E) = \alpha(E) /\gamma(E).$
Let $S\subset\R$ be a Borel set and assume that for all energies 
$E \in S,$ $\eta (E) \geq b > 0$. Then for any $g>0,$ there exists a constant $C_g$ 
such that for all $T>0$,
\begin{equation}\label{lb33}
   \langle \| e^{-iH_+ t} P_S \delta_1 \|^2_{C_g T^b}\rangle_T
\leq 1- \mu^+\left(\left\{E \left| \right.\eta(E) \geq b \right\} \right) +g. 
\end{equation}
\end{theorem}
\noindent{\it Remarks.} 1. As usual \cite{Guarneri,Last}, from the estimate \eqref{lb33}
follow lower bounds on the growth rate of moments of the position operator. \\
2. Theorem~\ref{pb} is a local version (and thus also a somewhat
stronger variant) of Theorem 1.2 of \cite{KL}. It can be proven directly by a slight
modification to the proofs in \cite{KL} (or see \cite{Tch2}), in which case it comes out
naturally as a multidimensional variant, where the exponent $\gamma(E)$ is connected
with decay of generalized eigenfunctions in the multidimensional space. We include this
theorem here, showing that it follows from Theorem~\ref{ldb}, mainly to illustrate the
connection of Theorem~\ref{ldb} with this type of results.\\

It is often convenient to filter scaling behaviors in terms of explicitly
defined scaling exponents and to formulate relations between scaling behaviors
as inequalities between such exponents. It is thus natural in our context
to define
\[
\overline\beta = \lim_{\delta\to 0}\;\limsup_{T\to\infty}
\frac{\log\left(\inf\left\{ L \left|\right. \langle \|e^{-iH_+ t} \delta_1 \|^2_{L} \rangle_T > \delta 
\right\}\right)}{\log T},
\]
\[
\underline\beta = \lim_{\delta\to 0}\;\liminf_{T\to\infty}
\frac{\log\left(\inf\left\{ L \left|\right. \langle \|e^{-iH_+ t} \delta_1 \|^2_{L} \rangle_T > \delta 
\right\}\right)}{\log T}.
\]
$\overline\beta$ and $\underline\beta$ are the upper and lower spreading rates associated
with the slowest spreading portion of the wavepacket. Local exponents for
the asymptotic scaling behavior of the scales $L^+_{\epsilon}(E)$ are given by
\[
\overline\lambda(E) = \limsup_{\epsilon\to 0}\frac{\log L^+_{\epsilon}(E)}{\log\epsilon^{-1}},\qquad
\underline\lambda(E) = \liminf_{\epsilon\to 0}\frac{\log L^+_{\epsilon}(E)}{\log\epsilon^{-1}}.
\]
Theorem \ref{hld} immediately implies
\[ 
\overline\beta\leq\mu^+\text{-ess}\inf\overline\lambda(E).
\] 
Moreover, we see that in problems where the solutions $u_0(n,E)$ behave nicely
enough to ensure that $\|u_0(n,E)\|_{L_1}/\|u_0(n,E)\|_{L_2}$ is small whenever
$L_1/L_2$ is small, Theorem \ref{ldb} would imply
\[ 
\mu^+\text{-ess}\inf\underline\lambda(E)\leq\underline\beta .
\] 
These inequalities are particularly interesting in cases where the problem exhibits
nicely scaling power law behaviors so that $u_0(n,E)$ behaves as described above and
also $\mu^+\text{-ess}\inf\underline\lambda(E)=\mu^+\text{-ess}\inf\overline\lambda(E)$.
If this happens, we see that we also have $\underline\beta = \overline\beta$ and the
spreading rate of the slowest spreading portion of the wavepacket is precisely given
by the slowest growth rate of the scales $L^+_{\epsilon}(E)$ with respect to the spectral
measure $\mu^+$.

Our second main goal in this paper is to apply the new dynamical upper bounds to 
the Fibonacci Hamiltonian, the most studied
of all one-dimensional models of quasicrystals. For this application we need a
whole-line (namely, $\ell^2(\Z)$) version of Theorem~\ref{hld}.
To formulate this whole-line version, 
note that the scales $L^-_\epsilon(E)$ and $\tilde{L}_{\epsilon}^- (E)$
can be defined in a way directly analogous to \eqref{ls}, \eqref{chartran}, involving 
the same kind of sums but taken over the negative half-line. Notice that the negative 
half-axis in our setting is $(\dots, -1,0),$ so that the summation in the analog of definition
\eqref{ls} will now start (or rather end) at $0$ instead of $-1.$ 
Let $\mu$ be the spectral measure of the whole-line operator $H$ corresponding to $\delta_1,$ defined 
by
\[ M(z)= \langle (H-z)^{-1} \delta_1, \delta_1 \rangle = \int_\reals \frac{d \mu(x)}{x-z}. \] 

\begin{theorem}\label{lindyn}
Let $H$ be a self-adjoint operator of the form \eqref{jacobi} on $\ell^2(\integers).$
Let the characteristic scales $L^{\pm}_{T^{-1}}(E)$ be defined by \eqref{ls}.
Then for any $T>0,$ $L_1,L_2>1,$ we have 
\begin{equation}\label{ub}
\langle \|e^{-iHt} \delta_1 \|^2_{L_1,L_2} \rangle_T > C \mu\left(\left\{ E \left| \right. 
L^{-}_{T^{-1}}(E) \leq L_1 \,\,\,{\rm and} \,\,\,L^{+}_{T^{-1}}(E) \leq L_2 \right\}\right),
\end{equation}
where $C$ is a universal positive constant. 
If $L_1,L_2>2$, the same result holds if $L^{\pm}_{T^{-1}}(E)$ is replaced in the statement
with $\tilde{L}^{\pm}_{T^{-1}}(E)$ defined by \eqref{chartran}.
\end{theorem}

The Fibonacci Hamiltonian $H_\lambda$ is defined by 
\[
(H_\lambda u)(n) = u(n+1) + u(n-1) + \lambda V(n) u(n),
\]
where $V(n) = \floor{(n+1) \omega } - \floor{ n\omega}$ and $\omega= (\sqrt{5}-1)/2$ 
is the golden mean. It is known \cite{BIST,Suto1,Suto2} that for every $\lambda > 0,$ 
$H_\lambda$ has purely singular continuous spectrum, and moreover, its spectrum (as a set)
is a Cantor set of zero Lebesgue measure. Lower bounds on wavepacket spreading rates for
$H_\lambda$ were
recently shown in \cite{JL2}.
We are going to show both upper and lower bounds
for the dynamics of $H_\lambda$ which imply that the spreading rate is intermediate
between ballistic ($\sim T$ at time $T$) and localized ($\sim T^0$). To the best of our
knowledge, this is the first proof of such behavior in an explicit model
of this type. 

\begin{theorem}\label{mainfib}
Let $H_\lambda$ be the Fibonacci Hamiltonian. Then 
\begin{enumerate}
\item[i)] There exists a constant $G>0$ such that for sufficiently large 
$\lambda,$ 
\begin{equation} \label{ubf}
\Big\langle
\|e^{-iH_\lambda t}\delta_1\|_{C (T+1)^{p_1(\lambda)},C (T+1)^{p_1(\lambda)}}^2 \Big\rangle_T
\geq  G \qquad \forall T>0, 
\end{equation}
where $p_1 (\lambda) = C_1 (\log \lambda)^{-1}(1+O((\lambda \log \lambda)^{-1})).$ 
\item[ii)] Given any $g>0,$ there exists a constant $C_g$ such that 
for every coupling $\lambda>0$ and every time $T>0,$ 
\begin{equation} \label{lbf}
\Big\langle \|e^{-iH_\lambda t}\delta_1\|_{C_g T^{p_2(\lambda)},C_g T^{p_2(\lambda)}}^2 \Big\rangle_T
\leq g, 
\end{equation}
where the positive function 
$p_2 (\lambda)$ satisfies $p_2(\lambda) = C_2 (\log \lambda)^{-1}(1+O( (\log \lambda)^{-1}))$ 
for large $\lambda$. 
\end{enumerate}
\end{theorem}
\noindent{\it Remarks.} 1. The first part of Theorem~\ref{mainfib} is an upper bound on wavepacket 
spreading for large coupling, which shows that on the average, there is a portion of the wavepacket remaining 
in a ball of the radius $\sim T^{p_1(\lambda)}$ at time $T.$ The second part of the 
Theorem 
provides a lower bound on wavepacket spreading, showing that on the average, only an arbitrarily
small part of the wavepacket remains in a ball of radius $\sim T^{p_2(\lambda)}$ at time $T.$
The constants $p_1(\lambda)$ and $p_2(\lambda)$ have the same (logarithmic in $\lambda$) 
asymptotic behavior for large coupling, up to a constant factor in front of the main term. 
This implies that the logarithmic law is precise and cannot be improved in the estimates
\eqref{ubf}, \eqref{lbf}. \\
2. It will follow from the proof below that $\lambda >8$ is sufficient to get a nontrivial upper 
bound ($p_1(\lambda)<1$ in \eqref{ubf}). This range is not optimal and can be improved by 
additional technical effort; however our current methods do not allow to extend the bound all the way 
to $\lambda=0.$ Moreover, if $\lambda$ is so small that $p_1(\lambda)\geq 1$, then there is little interest
since this would constitute a ballistic (or worse) bound. \\ 
3. We remark that Sinai, in a recent paper \cite{Sinai}, studied anomalous transport
(in terms of moments of the position operator) for an almost periodic potential in a different
setting (in fact, in a pure point regime). \\

Although we only treat discrete operators in this paper, results parallel
to Theorems~\ref{hld}, \ref{hldt}, \ref{ldb}, \ref{pb}, \ref{lindyn} also hold in continuous
settings by direct translation of the arguments given here.

The rest of this paper is organized as follows. In Section 2 we prove some auxiliary facts relating the 
scale $L^+_\epsilon(E)$ and the $m$ function.
In Section 3 we derive the
upper bounds on wavepacket spreading. In particular, the proofs of
 Theorem~\ref{hld}, Theorem~\ref{lindyn} and Corollary~\ref{hldt} appear
 there.
In Section 4
we show the new lower bounds, Theorems~\ref{ldb} and \ref{pb}.
In Sections 5 and 6 we treat the Fibonacci Hamiltonian and prove Theorem~\ref{mainfib}. Finally,
In Section 7 we discuss how our results are applicable to multidimensional problems.

\section{Bounds on the $m$ function}  

We are going to prove a series of auxiliary estimates relating the behavior of solutions 
$u_{0,\pi/2}$ and the half-line $m$ functions. The line of the argument follows Jitomirskaya-Last 
extension of subordinacy theory \cite{JL1}, but there will be an essential modification 
that will be crucial for the derivation of the dynamical criteria.

Recall the following simple example of variation of parameters which will prove extremely useful:

\begin{lemma}\label{vplem}
Suppose $w(n,z),w(n,E)$ solve \eqref{ef} for spectral parameters $z,E$ respectively.
Assume that $w(0,z)=w(0,E)$ and $w(1,z)=w(1,E).$ Then for $n \geq 0,$ 
\begin{equation}
w(n,z) = w(n,E) + (z-E) \sum_{m=1}^n k(n,m;E) w(m,z)   \label{vpeq}
\end{equation}
where
$$
k(n,m;E)=u_0(n,E)u_{\pi/2}(m,E)-u_{\pi/2}(n,E)u_0(m,E)
$$
and with the convention that the sum is zero if the upper limit is less than the lower one. 
\end{lemma}

\begin{proof}
It is a direct computation to check that the right-hand side satisfies \eqref{ef} with
the spectral parameter $z$.
Hence, on both sides we have solutions of \eqref{ef} which coincide at the sites 
$0,1.$ This implies the equality. 
\end{proof}

Equation \eqref{vpeq} has the form
\[
w(n,z) = w(n,E)+(z-E) (K(E) w)(n,z),
\]
where $K(E)$ is an integral operator with the kernel $k(n,m;E)$.
As introduced in \eqref{normL}, $\|\cdot\|_L$ 
defines a norm on an $\floor{L}+1$-dimensional Hilbert space
(except when $L$ is an integer, in which case the dimension is $L$).  
The corresponding inner product is
\begin{equation}\label{inproL}
 \langle \phi_1, \phi_2 \rangle_L = \sum\limits_{n=1}^{\floor{L}} \phi_1(n)\phi_2(n) +
(L-\floor{L})\phi_1(\floor{L}+1)\phi_2(\floor{L}+1).  
\end{equation}
We wish to estimate the norm of the operator $K(E)$ acting in this space:

\begin{lemma}\label{twosc} The Hilbert-Schmidt norm of $K(E)$ is given by
\begin{align} 
|||K(E)|||^2_L  &= \| u_0 \|^2_L \|u_{\pi/2}\|^2_L - | \langle u_0, u_{\pi/2} \rangle_L |^2 \label{Kb} \\
 \label{2sceq}  &= \max_\theta \|u_\theta (n,E)\|^2_{L} \,\cdot\, \min_\theta \|u_\theta (n,E)\|^2_{L}.
\end{align}
where $u_0,u_{\pi/2}$ are the solutions at energy $E$. In particular,
this gives an upper bound on the operator norm of $K$.
\end{lemma}

\begin{proof}
Recall that the Hilbert-Schmidt norm of an integral operator is equal to the $L^2$ norm of its kernel.
For integer $L$,
$$
|||K(E)|||^2_L = \sum_{n,m=1}^{L} |k(n,m;E)|^2 = \| u_0 \|^2_L \|u_{\pi/2}\|^2_L - | \langle u_0, u_{\pi/2} \rangle_L |^2
$$
can be shown fairly easily. Non-integer $L$ merely require a more lengthy computation.

To show the second equality, consider the $2\times 2$ matrix
$$
Q = \begin{bmatrix} \| u_0 \|^2_L & \langle u_0, u_{\pi/2} \rangle_L \\
        \langle u_{\pi/2}, u_0 \rangle_L & \| u_{\pi/2} \|^2_L \end{bmatrix}.
$$
It is self-adjoint and positive, and if $\vec{e}_\theta=(\cos \theta, \sin \theta)$, then
$$
\vec{e}_\theta Q {\vec{e}_\theta}^{\;t} = \|u_\theta\|_L^2.
$$
The right-hand side of \eqref{Kb} is the determinant of $Q$ and so the product of the eigenvalues of $Q$.
By the minimax characterization, these eigenvalues are given by the factors in \eqref{2sceq}.
\end{proof}

\it Remark. \rm Notice that if the maximum of $\|u_\theta(n,E)\|_L$ is achieved at 
$\theta_0,$ the minimum is achieved at the orthogonal boundary condition $\theta_0^\bot$
($u_{\theta_0^\bot}$ satisfies $u_{\theta_0^\bot}(0,E)= \cos \theta_0,$ $u_{\theta_0^\bot}(1,E)=
-\sin \theta_0$). This is because the eigenvectors of the self-adjoint matrix $Q$ are orthogonal.

{From} \eqref{2sceq} we know that $|||K(E)|||_L$ is strictly increasing.  Hence
$\epsilon |||K(E)|||_L = 1$ determines $L$ as a function of $\epsilon$ and $E$.  As in the introduction,
we denote this length scale by $\Lep(E)$.  The main result of this section is an estimate for
the norm of the Weyl solution $u_+$ over this length scale. Before stating and proving this result,
let us recall some facts about the Weyl theory of Jacobi matrices (see, e.g., \cite{Ber}):
For each $z\in\complex\setminus\reals$, there exists a unique solution 
of \eqref{ef} which is square summable and obeys $u_+(0,z)=1$.  One may write this
solution in the form
$$
u_+(n,z) = u_{\pi/2}(n,z)- m_+(z) u_0(n,z),
$$
and so define the Weyl $m$-function $m_+(z)$.  As we will discuss a little more in the next section, $m_+(z)$
captures all of the spectral information about $H_+$.  For the moment, however, we merely need the
observation that 
\begin{equation}\label{sp}
a(0)\Im m_+(z) = \Im z \sum_{n=1}^\infty |u_+(n,z)|^2,
\end{equation}
which follows from summation by parts.
Recall that we assume $a(0)=1$.

\begin{theorem}\label{mainm}
Fix $E \in \reals$ and $\epsilon >0$, then we have
\begin{equation}\label{jlb}
2-\sqrt{3}  \leq \frac{\|u_0\|_\Lep |m_+(E+i\epsilon)|}{\|u_{\pi/2}\|_\Lep} \leq  2+\sqrt{3}.
\end{equation}
Moreover, 
\begin{equation}\label{betam}
\epsilon \|u_+\|^2_\Lep \geq \frac{|m_+(E+i\epsilon)|}{  4\epsilon \|u_0\|_\Lep \|u_{\pi/2}\|_\Lep  }
        \geq \tfrac{2-\sqrt{3}}{16} \Im m_+(E+i\epsilon).
\end{equation}
In both formulae, $u_0,u_{\pi/2}$ are the solutions at energy $E$.
\end{theorem}

\noindent {\it Remark}. The bounds \eqref{jlb} are very similar to those proved by
Jitomirskaya-Last \cite{JL1}. The difference is that the scale $L_\epsilon^+$ is given
in \cite{JL1} by the condition 
\[
\epsilon \|u_0\|_{L^+_\epsilon}\|u_{\pi/2}\|_{L^+_\epsilon} = 1.
\]
The scale $\Lep$ which we define by
\begin{equation}\label{scale}
|||K|||_{\Lep} = \max_\theta \|u_\theta (n,E)\|^2_{\Lep} \,\cdot\, \min_\theta \|u_\theta (n,E)\|^2_{\Lep}
                = \epsilon^{-1}
\end{equation}
is larger (compare \eqref{Kb}), and the constants appearing in \eqref{jlb}
are sharper than those in \cite{JL1}. The scale defined by \eqref{scale} might be less 
convenient for dimensional spectral analysis, since its definition is more involved.
However, an improvement contained in \eqref{scale} will be quite crucial for the proof 
of dynamical criteria in the next section.

\begin{proof} To shorten the formulae which will follow, we introduce
\begin{equation*}
\beta = \frac{ \langle u_0, u_{\pi/2} \rangle_\Lep }{\|u_0\|_\Lep \|u_{\pi/2}\|_\Lep}
        \quad \text{ and } \quad \zeta=\sqrt{1-\beta^2}.
\end{equation*}
Notice that Cauchy-Schwarz says $\zeta>0$. At times we shall also write $m_+$ for $m_+(E+i\epsilon)$.
Notice that by Lemma~\ref{twosc}, \eqref{Kb} and \eqref{ls}, the scale 
$L_\epsilon^+ (E)$ given by \eqref{ls} is chosen exactly in a way to ensure that
$\|K(E)\|_L \leq \epsilon^{-1}$ if $L \leq L_\epsilon^+ (E).$
Lemmas~\ref{vplem},\ref{twosc} and the definition of $\Lep$ combine to show that
\begin{equation}\label{startm}
 4\|u_+(n,E+i\epsilon)\|^2_{L^+_\epsilon} \geq \|u_{\pi/2}(n,E) - m_+(E+i\epsilon) u_0(n,E)\|^2_{L^+_\epsilon}. 
\end{equation}
Using \eqref{sp}, we see that 
$$
4\epsilon^{-1} \Im m_+ \geq 4\|u_+\|^2_\Lep \geq \|u_{\pi/2}\|^2_{L^+_\epsilon} + |m_+|^2\cdot\|u_0\|^2_{L^+_\epsilon}
        - 2 \Re m_+ \langle u_0, u_{\pi/2} \rangle_{L^+_\epsilon}.
$$
The definitions of $\beta,\zeta$ and $\Lep$ are such that $\epsilon\zeta\|u_0\|_{\Lep}\|u_{\pi/2}\|_{\Lep}=1$.
Thus, 
\begin{equation}\label{keym}
\Im m_+ \geq  \epsilon\|u_+\|^2_{L^+_\epsilon}
 \geq \frac{1}{4\zeta} \bigg\{ \frac{\|u_{\pi/2}\|_\Lep}{\|u_0\|_\Lep} + 
|m_+(E+i\epsilon)|^2 \frac{\|u_0\|_\Lep}{\|u_{\pi/2}\|_\Lep} - 2\beta\Re m_+ \bigg\}.
\end{equation}
This equation implies both claims of the theorem. To prove \eqref{jlb}, notice that $|\pm\beta+i\zeta|=1$  and
so
$$
-2|\beta \Re(m_+)| - 4\zeta\Im m_+ \geq -2|m_+| - 2\Re\big[(\pm\beta+i\zeta)m_+\big] \geq -4|m_+|,
$$
where the sign is chosen to make $\pm\beta\Re(m_+)\geq 0$.  Applying this to \eqref{keym}, leads to
$$
\frac{\|u_0\|_\Lep}{\|u_{\pi/2}\|_\Lep} |m_+(E+i\epsilon)|^2 - 4|m_+(E+i\epsilon)| +
        \frac{\|u_{\pi/2}\|_\Lep}{\|u_0\|_\Lep} \leq 0,
$$
which is equivalent to \eqref{jlb}. To prove \eqref{betam}, use \eqref{keym} to obtain
\begin{align*}
\epsilon \|u_+(n,E+i\epsilon)\|^2_{L^+_\epsilon}
 &\geq \frac{|m_+|}{4\zeta} \bigg\{ \frac{\|u_{\pi/2}\|_\Lep}{|m_+|\|u_0\|_\Lep}+
        \frac{\|u_0\|_\Lep |m_+|}{\|u_{\pi/2}\|_\Lep} - 2\beta \frac{\Re m_+}{|m_+|} \bigg\}    \\
&\geq \frac{|m_+|}{4\zeta} (2-2|\beta|) = \frac{\zeta|m_+|}{2(1+|\beta|)}       \\
&\geq \tfrac{1}{4} \zeta|m_+(E+i\epsilon)|.
\end{align*}
(We used $x+x^{-1} \geq 2$ in the second step.)
This proves the left-hand inequality in \eqref{betam}, because
$\epsilon\zeta\|u_0\|_{\Lep}\|u_{\pi/2}\|_{\Lep}=1$.  To complete the proof,
we need to show the right-hand inequality.
{From} \eqref{startm} and \eqref{sp} we infer that 
$$
4\Im m_+ \geq \epsilon (\Im m_+)^2 \|u_0\|^2_{L^+_\epsilon},
$$
and from \eqref{jlb}, 
$$
\|u_0\|^2_{L^+_\epsilon} \geq
(2-\sqrt{3}) |m_+(E+i\epsilon)|^{-1}  \|u_0\|_{L^+_\epsilon}\|u_{\pi/2}\|_{L^+_\epsilon}.
$$
Combining these two gives the right-hand side of \eqref{betam}.
\end{proof}

The following variant of the inequality 
\eqref{betam} in Theorem~\ref{mainm} will prove useful in the next section.

\begin{proposition}\label{BigTwo}
Suppose $E,E'\in\reals$, $\epsilon>0$ and $|E-E'|<\epsilon.$
Then there exists a universal
constant $C$, so that
\begin{equation}\label{Cbound}
\epsilon \|u_+(n,E'+i\epsilon)\|^2_{L^+_\epsilon(E)} \geq C \Im m_+(E'+i\epsilon).
\end{equation}
\end{proposition} 

\noindent \it Remark. \rm It is important that while the bound in \eqref{Cbound} involves
$m_+$ and the solution
$u_+$ at the energy $E'+i\epsilon,$ the scale $L^+_\epsilon$ is defined at the energy $E.$

\begin{proof} This result is a direct corollary of the proof of Theorem~\ref{mainm}. 
All statements of Theorem~\ref{mainm} hold with $E+i\epsilon$ replaced by $E'+i\epsilon$ and adjusted
constants. Specifically,
the only change introduced by replacing $u_+(n,E+i\epsilon)$ by $u_+(n,E'+i\epsilon)$, $|E-E'|<\epsilon$,
is a change in the constant in \eqref{startm} from $4$ to $(1+\sqrt{2})^2 \leq 6$.
Following through the proof with the constant $6$ shows that $C$ may be chosen to be
$(3-2\sqrt{2})/36.$
\end{proof}

\section{An upper bound on wavepacket spreading}

We begin by relating the dynamical quantity we need to estimate to a solution of equation \eqref{ef}.
As before,  
\[ u_+(n,z)= u_{\pi/2}(n,z) - m_+(z)u_{0}(n, z) \]
denotes the unique $\ell^2(\integers_+)$ solution for $z \in \complex \setminus \reals$ with $u_+(0,z)=1$. 


\begin{lemma}\label{rese}
For every $z \in \complex \setminus \reals$ and $n \geq 1,$ 
\begin{equation}\label{resed}
 \langle(H_+-z)^{-1}\delta_1, \delta_n \rangle = - u_+(n,z). 
\end{equation}
\end{lemma}
\begin{proof}
It's easy to see from the definition of the resolvent of $H_+$ that
$\langle(H_+-z)^{-1}\delta_1, \delta_n \rangle$, $n\geq 1$, is an $\ell^2(\integers_+)$
solution of \eqref{ef}. It thus follow from the uniqueness of the Weyl solution that
it must be a multiple of $u_+(n,z)$. The Lemma is thus implied by the
definition of $m_+(z)$ and our normalization of $u_+(n,z)$.
\end{proof}

\begin{lemma}\label{dynsol}
For any $T>0,$ 
\begin{equation}\label{dynsoleq}
\langle |\langle e^{-iH_+ t} \delta_1, \delta_n \rangle|^2 \rangle_T = \frac{1}{\pi T} \int\limits_\reals 
|u_+(n, E+ \frac{i}{T})|^2 \, dE. 
\end{equation}
\end{lemma}
\begin{proof}
Recall that 
\[ \langle |\langle e^{-iH_+ t} \delta_1, \delta_n \rangle|^2 \rangle_T =
\frac{2}{T} \int\limits_0^\infty e^{-2t/T} | \langle e^{-iH_+ t} \delta_1, \delta_n \rangle |^2\, dt. \]
The right-hand side in the above equality is a constant times the square of the  $L^2$ norm of the function 
$e^{-t/T} \int\limits  e^{-iE' t}\, d\mu_{1,n} (E'),$ where the complex measure $\mu_{1,n}$ is defined by 
\[ \langle (H_+ -z)^{-1} \delta_1, \delta_n \rangle = \int_\reals \frac{d \mu_{1,n}(x)}{x-z}. \]
The Fourier transform of this function is equal to 
\[ \frac{i}{2\pi} \int\limits_\reals \frac{d \mu_{1,n}(E')}{(E-E')+(i/T)}= -\frac{i}{2\pi}\langle
(H_+ -E-\frac{i}{T})^{-1}\delta_1,
\delta_n \rangle = \frac{i}{2\pi}u_+(n,E+\frac{i}{T}) \]
by \eqref{resed}.
Now \eqref{dynsoleq} follows from Parseval's equality for the Fourier transform.
\end{proof}
Lemma~\ref{dynsol} implies, 
\begin{equation}\label{pars}
\langle \| e^{-iH_+ t} \delta_1 \|^2_L \rangle_T = \frac{1}{\pi T} \int\limits_\reals \|u_+(n,E+\frac{i}{T})\|^2_L \, dE.
\end{equation}

Given a set $S \subset \reals,$ let us denote by $S_\epsilon$ the $\epsilon$-neighborhood of the set $S.$ 
We now complete the proof of the half-line dynamical bound.
\begin{proof}[Proof of Theorem~\ref{hld}]
Let 
\[ S = \left\{ E \left| \right. L^+_{T^{-1}}(E)\leq L \right\}. \]
By Lemma~\ref{dynsol}, to 
prove Theorem~\ref{hld} it suffices to bound the right-hand side of 
\eqref{pars} from below by $C\mu^+(S).$ Hence, it suffices to estimate from below
the integral 
\[ \epsilon \int\limits_{S_\epsilon} \|u_+(n,E+i\epsilon)\|_L^2 \, dE, \]
where $L \geq L^+_\epsilon(E)$ for $E \in S$ (we think of $\epsilon=T^{-1}$). 
By the definition of $S_\epsilon,$ for every $E' \in S_\epsilon$ there is $E \in S$ 
such that $|E -E'| <\epsilon.$ Applying Proposition~\ref{BigTwo}, we obtain
\begin{equation} \label{BT}
\epsilon \int\limits_{S_\epsilon} \|u_+(n,E+i\epsilon)\|_L^2 \, dE \geq 
C\int\limits_{S_\epsilon} \Im m_+(E+i\epsilon) \, dE. 
\end{equation}
By the Fubini theorem, the left-hand side in \eqref{BT} is equal to 
\[ C\int\limits_\reals d\mu_+(x) \int\limits_{S_\epsilon} \frac{\epsilon \, dE}{(x-E)^2+\epsilon^2}
\geq C\int\limits_S d \mu_+(x) \int\limits_{-\epsilon}^{\epsilon} \frac{\epsilon \,dE'}{(E')^2 +
\epsilon^2} \geq C \mu_+(S). \]
\end{proof}

The following Lemma shows that Theorem~\ref{hld} implies Corollary~\ref{hldt}.

\begin{lemma}\label{trans}
For every $E \in \R$ and $L\geq 2$, 
\begin{equation}\label{transeq}
4 \|\Phi (1,E)^{-1}\|^2 \left( \max_\theta \|u_\theta (n,E)\|^2_{L} \right)
\left( \min_\theta \|u_\theta (n,E)\|^2_{L} \right) \geq \|\Phi(E)\|_L^2.  
\end{equation}
\end{lemma}
\begin{proof}
A direct computation using the definition of $\Phi(n,E)$ shows that for any $\theta,$ 
\[ \|\Phi (n,E)\|^2 \leq (u_\theta (n,E))^2+(u_\theta (n+1,E))^2+(u_{\theta^\bot} (n,E))^2+
(u_{\theta^\bot} (n+1,E))^2. \]
Thus,
\[
\|\Phi(E)\|_L^2 \leq 2(\|u_\theta(n,E)\|^2_L+\|u_{\theta^\bot}(n,E)\|^2_L)
\leq 4 \left( \max_\theta \|u_\theta (n,E)\|^2_{L} \right).
\]
For $L\geq 2$, we have for any $\theta$,
$\|u_\theta\|_L\geq \|u_\theta\|_2\geq \|\Phi (1,E)^{-1}\|^{-1}$, and we thus see
that \eqref{transeq} holds.
\end{proof}

\begin{proof}[Proof of Corollary~\ref{hldt}]
Since $\tilde{L}_\epsilon^+(E)$ was defined in \eqref{chartran} by
\[ \|\Phi(E)\|_{\tilde{L}_{\epsilon}^+ (E)} = 2 \|\Phi (1,E)^{-1}\| \epsilon^{-1}, \]
we see from \eqref{ls} and Lemma~\ref{trans}, that
$\tilde{L}_\epsilon^+(E) \geq L_\epsilon^+(E)$ whenever $L_\epsilon^+(E)\geq 2$.
Therefore, Corollary~\ref{hldt} follows immediately from Theorem~\ref{hld}.
\end{proof}


It remains to consider the whole line case. 
Denote by $u_-(n,z),$ $z \in \complex \setminus \reals$,
a solution of \eqref{ef} which belongs to $\ell^2(\integers_-),$ $Z_-=\{0,-1,\dots,\},$ 
and is normalized by 
\[ u_-(n,z) = u_{\pi/2}(n,z) +m_-(z)u_0(n,z). \]
Recall the notation $H$ for the operator defined by \eqref{jacobi} on
$\ell^2(\integers).$ 

\begin{proof}[Proof of Theorem~\ref{lindyn}]
As in Lemma~\ref{dynsol}, we can show 
$$
\langle \|e^{-iHt} \delta_1 \|_{L_1,L_2}^2 \rangle_T = \frac{1}{\pi T} \int_\reals
\|G(1,n,E+\frac{i}{T})\|^2_{L_1,L_2}dE,
$$
where 
\[ G(1,n,z) = \langle (H-z)^{-1} \delta_1, \delta_n \rangle. \]
One can easily verify that
\[ G(1,n,z) = \left\{ \begin{array}{ll} \frac{m_-(z)}{a(0) (m_+(z)+m_-(z))} u_+(n,z) & n \geq 1 \\
 \frac{-m_+(z)}{a(0) (m_+(z)+m_-(z))} u_-(n,z) & n < 1. \end{array} \right. \]
In particular, 
\[ G(1,1,z)= M(z)=\frac{- m_+(z) m_-(z)}{a(0) (m_+(z) + m_-(z))}, \]
and the spectral measure $\mu$ of $H$ corresponding to the vector $\delta_1$ satisfies
\[ M(z) = \int\limits_\reals \frac{d \mu(x)}{x-z}. \]
Now define
\[ S =  \left\{ E \left| \right. L^+_{T^{-1}}(E)\leq L_2 \,\,\,{\rm and} \,\,\, 
L^-_{T^{-1}}(E)\leq L_1  \right\}. \]
Repeating the same arguments that led us to the proof of Theorem~\ref{hld}, we get 
($\epsilon = T^{-1}$) 
\[ \epsilon \|G(1,n,E+i\epsilon)\|_{-1,L_2}^2 \geq C \int\limits_{S_\epsilon} \frac{|m_-(E+i\epsilon)|^2}
{|m_+(E+i\epsilon)+m_-(E+i\epsilon)|^2} \Im m_+(E+i\epsilon) \, dE \]
and similarly
\[ \epsilon \|G(1,n,E+i\epsilon)\|_{L_1,0}^2 \geq C \int\limits_{S_\epsilon} \frac{|m_+(E+i\epsilon)|^2}
{|m_+(E+i\epsilon)+m_-(E+i\epsilon)|^2} \Im m_-(E+i\epsilon) \, dE. \]
Combining these two inequalities, we get 
\[ \epsilon \|G(1,n,E+i\epsilon)\|_{L_1,L_2}^2 \geq C \int\limits_{S_\epsilon} \Im M(E+i\epsilon)\,dE \geq
C \mu(S), \]
similarly to the last step in the proof of Theorem~\ref{hld}.
\end{proof}

\section{A lower bound on wavepacket spreading}

Before treating our main application, the Fibonacci Hamiltonian, we pause to 
prove the general lower bound on dynamics given by Theorem~\ref{ldb}.
\begin{lemma}\label{evol}
Let $P_S$ be a spectral projection of the operator $H_+$ on a measurable set $S.$ Then 
\begin{equation}\label{evoleq}
\left( e^{-iH_+ t} P_S \delta_1 \right) (n)= \int\limits_S e^{-iEt} u_0 (n,E) d\mu^+(E).
\end{equation}
\end{lemma}
\begin{proof}
Since $\delta_1$ is a cyclic vector for $H_+$, it follows from the spectral theorem that $H_+$
is unitarily equivalent to multiplication by the parameter on $L^2(\R,d\mu^+)$. $u_0 (n,E)$, $n\geq 1$,
is known (see, e.g., \cite{Ber}) to be the representation of $\delta_n$ in this space.
\end{proof}

\begin{proof}[Proof of Theorem~\ref{ldb}]
Notice that according to Lemma~\ref{evol},
\begin{align}
& \big\langle \| e^{-iH_+ t } P_S \delta_1 \|_L^2 \big\rangle_T
   \nonumber \\
\leq {}& \frac{2}{T} \int\limits_0^\infty e^{-2t/T} \sum\limits_{n \leq L}
\int\limits_S
e^{iEt} u_0 (n,E)\, d\mu^+(E)  \int\limits_S
e^{-iE't} u_0 (n,E')\, d\mu^+(E') \, dt \nonumber  \\
= {}& \sum\limits_{n \leq L} \int\limits_S  \int\limits_S
u_0 (n,E')u_0(n,E) \frac{4T^{-2}}{(E-E')^2 + 4T^{-2}}
   \,d\mu^+(E)\,d\mu^+(E') \nonumber \\
\label{oyp}
\leq {}& \int\limits_S  \left( \int\limits_S  \frac{4T^{-2}}{(E-E')^2 +4T^{-2}}
   \,d\mu^+(E') \right)
\|u_0 (n,E)\|_L^2 \,d\mu^+(E).
\end{align}
We used Cauchy-Schwartz in the last step with respect to the product measure $d\mu^+(E) d\mu^+(E').$ 
Notice that the term in the brackets in \eqref{oyp} does not exceed (setting $\epsilon=T^{-1}$)
\[ 2\epsilon \Im m_+(E+2i\epsilon) \leq 4\epsilon \Im m_+(E+i\epsilon). \]
By the estimates \eqref{jlb} and \eqref{betam} of Theorem~\ref{mainm}, we have 
\begin{equation}\label{mmest}
\epsilon \Im m_+(E+i\epsilon) \|u_0(n,E)\|^2_{L^+_\epsilon} \leq (2+\sqrt{3}) \frac{\epsilon \Im m_+(E+i\epsilon)
\|u_0\|_{L^+_\epsilon}\|u_{\pi/2}\|_{L^+_\epsilon}}{|m_+(E+i\epsilon)|} \leq C.
\end{equation}
Combining \eqref{oyp} and \eqref{mmest}, we obtain 
\[  \langle \| e^{-iH_+ t } P_S \delta_1 \|_L^2 \rangle_T \leq C \int\limits_S \frac{ \|u_0 (n,E)\|_L^2}
{ \|u_0(n,E)\|^2_{L^+_\epsilon}} d \mu^+(E), \]
which is exactly what we wanted to show.
\end{proof}

\begin{proof}[Proof of Theorem~\ref{pb}]
By the extension of subordinacy theory due to Jitomirskaya and Last \cite{JL1},
\[ \liminf_{L \rightarrow \infty} \frac{\|u_0\|_L^{2-\alpha}}{\|u_{\pi/2}\|_L^\alpha}=0 
\Leftrightarrow D^\alpha \mu^+(E) = \infty. \]
Therefore, for a.e.~$E,$ we have with some $C_1(E)>0$
\begin{equation}\label{21b}
 \|u_0\|_L^{\frac{2-\alpha(E)}{\alpha(E)}} \geq C_1(E) \|u_{\pi/2}\|_L. 
\end{equation}
By \eqref{21b}, the definition \eqref{ls} of $L^+_\epsilon$ and relation \eqref{Kb}, 
\[ \epsilon^{-1} \leq 
\|u_0\|_{L^+_\epsilon}\|u_{\pi/2}\|_{L^+_\epsilon} \leq  C_1^{-1}(E)\|u_0\|_{L^+_\epsilon}^{\frac{2}{\alpha(E)}}. \]
Therefore, 
\begin{equation}\label{plb}
\|u_0\|_{L^+_\epsilon}^2 \geq C_1^{-1}(E) \epsilon^{-\alpha(E)}.
\end{equation}
On the other hand, by the definition of $\gamma(E),$ we have 
\[ \|u_0\|_L^2 \leq C_2(E) L^{\gamma(E)}. \]
Denote 
\[ S = \left\{ E \left| \right. \eta(E) \geq b \right\} \]
 (recall $\eta(E)=\alpha(E)/\gamma(E)$).
Let $S_1 \subset S$ be the set such that for $E \in S_1,$ 
$C_1(E)C_2(E) \leq C$
with some uniform constant $C.$ 
Clearly, $S_1$ may be chosen so that 
 $\| P_{S \setminus S_1} \delta_1\|^2 = \mu^+(S \setminus S_1)$ 
 is as small as we want by adjusting $C.$    
Then from \eqref{eqldb} of Theorem~\ref{ldb} we find (with $C$ denoting different universal constants
in different places) 
\begin{align}\label{vls}
   \langle \| e^{-iH_+ t } P_{S_1} \delta_1 \|_L^2 \rangle_T &\leq C
\int\limits_{S_1} \frac{ \|u_0 (n,E)\|_L^2}
{ \|u_0(n,E)\|^2_{L^+_\epsilon}} d \mu^+(E) \\ \nonumber
&\leq C  \int\limits_{S_1} T^{-\alpha(E)} L^{\gamma(E)} d\mu^+(E). 
\end{align}
Choosing $L=C_g T^b,$ we can make the left hand side in \eqref{vls} arbitrarily small by changing $C_g.$ 
It remains to observe that 
\begin{align*}
\langle \| e^{-iH_+ t } \delta_1 \|_L^2 \rangle_T &\leq 
\langle (\| e^{-iH_+ t } P_{S_1} \delta_1 \|_L+\| e^{-iH_+ t }(I-P_{S_1})
   \delta_1 \|_L)^2 \rangle_T \\
&\leq \left( \sqrt{\langle \| e^{-iH_+ t } P_{S_1} \delta_1 \|_L^2
   \rangle_T} + \|(I-P_{S_1})\delta_1\| \right)^2.
\end{align*}
We claim that by choosing $S_1$ and then $C_g$ the last expression can be made smaller than $1-\mu^+(S)+g$ 
for any $g>0.$ Indeed, first choose $S_1$ so that 
\[  \|(I-P_{S_1})\delta_1\|^2 < 1-\mu^+(S) +(g/2). \]
 Then choose $C_g$ so that 
\[ \langle \| e^{-iH_+ t } P_{S_1} \delta_1 \|_L^2 \rangle_T \leq (g/5) \]
(assuming $g$ is small). This completes the proof, demonstrating that the part of the wavepacket 
corresponding to the energies in $S$ leaves (on the average) a ball of the size $\sim T^b$ at time $T.$  
\end{proof}

\section{The Fibonacci Hamiltonian: bounds on traces}

Let us recall from the Introduction that the Fibonacci Hamiltonian is the discrete
Schr\"odinger operator 
\begin{equation}
  [H_\lambda u](n) = u(n+1) + u(n-1) + \lambda V(n)u(n)
\end{equation}
acting on $\ell^2(\integers)$ with the Fibonacci potential $V(n)=\floor{(n+1)\omega}-\floor{n\omega}$.
Here $\omega$ is the golden ratio, $(\sqrt{5}-1)/2.$
The most important property of this potential is the substitution rule \cite{Suto1}, 
\begin{equation}
\label{subs}
V(q_{k}+n) = V(n) \quad \text{ for $n=1,...,q_{k}$ and  $k\geq 3$}
\end{equation}
where $q_k$ denote the Fibonacci numbers, $q_0=1,q_1=1,q_k=q_{k-1}+q_{k-2}$.

\noindent
For $k \geq 1$ we define by $\Phi_k$ the transfer matrix $\Phi(q_k,E):$ 
$$
\begin{bmatrix} u(q_k+1) \\ u(q_k) \end{bmatrix} = \Phi_k \begin{bmatrix} u(1) \\ u(0) \end{bmatrix}
$$
for all $u(n)$ satisfying $H_\lambda u=Eu$. It is convenient to make the additional definitions
\begin{equation*}
\Phi_{-1} = \begin{bmatrix} 1 & -\lambda \\ 0 & 1 \end{bmatrix}, \quad
\Phi_{ 0} = \begin{bmatrix} E &     -1   \\ 1 & 0 \end{bmatrix}.
\end{equation*}
{From} the substitution rule it follows that $\Phi_{k+1}=\Phi_{k-1} \Phi_k,$ $k \geq 0$ \cite{KKT}. 
Let us denote $x_k(E) = \tr \Phi_k.$  
With more work, one
obtains the trace map and trace invariant \cite{KKT}:
\begin{gather}
\label{TM}
x_{k+1} = x_{k} x_{k-1} - x_{k-2} \\
\label{Invar}
x_{k+1}^2  + x_{k}^2 + x_{k-1}^2 - x_{k+1}x_{k}x_{k-1}= 4+\lambda^2.
\end{gather}
We denote by $\sigma_k$ the spectrum of the periodic potential with period $q_k$
coinciding with $V(n)$ for $n=1, \dots, q_k.$ By the Bloch theorem, $\sigma_k$ is a set of 
intervals (bands) for which $x_k(E) \in [-2,2]$ (see, e.g. \cite{Toda}). 
Moreover, in each band $x_k(E)$ varies monotonically in $[-2,2]$ and takes values 
$\pm 2$ at the ends.  
Due to the relations \eqref{TM}, \eqref{Invar}, 
the traces are among the most convenient objects of study 
in the Fibonacci model. Since $x_k$ describe the spectrum of periodic approximants,
there is also a natural relation to the spectrum of the limiting Fibonacci Hamiltonian.
The following properties of the traces and their relation to the spectrum of $H_\lambda$ 
are quite useful:
\begin{proposition}\label{background}
\begin{enumerate}

\item[i)] The spectrum of $H_\lambda$ coincides with the set of energies for which 
the sequence $x_k(E)$ is bounded.

\item[ii)] If $|x_k(E)|>2,$ $x_{k+1}(E)>2$ for some $k,$ 
then the sequence $x_n(E)$ is unbounded. 

\item[iii)] If $\lambda >4,$ there cannot exist $E,k$ such that 
$|x_k(E)| \leq 2$, $|x_{k+1}(E)| \leq 2$ and $|x_{k+2}(E)| \leq 2.$

\end{enumerate}
\end{proposition} 
\begin{proof}
The first and second statements have been proved by S\"ut\H{o} \cite{Suto1}. 
The third statement is a direct consequence of the trace invariant \eqref{Invar}.
\end{proof}

Our goal in this section is to prove a lower bound on the derivative $|x'_k(E)|$
for energies in the spectrum of $H_\lambda.$  
In the next section we show that this bound translates directly into the lower bound 
on the growth of $\|\Phi(n,E)\|_L^2$ that we need for the upper bounds on dynamics. 
The results of this section are contained, in a somewhat different form,
 in the preprint of Raymond \cite{Ray},
who used them to derive an upper bound on the Hausdorff dimension of the spectrum
(which alone does not imply any upper bounds on dynamics, see \cite{KL}). 
For the sake of completeness, we present here a simplified version of the argument 
given in \cite{Ray}. The main result of this section is 

\begin{proposition}\label{tracebound}
Assume that the coupling $\lambda$ is sufficiently large ($\lambda >8$ will do).
Then for every $E$ in the spectrum of $H_\lambda,$ the derivative of the trace 
$x_k(E)$ satisfies 
\begin{equation}\label{tracelb}
|x'_k(E)| \geq \xi (\lambda)^{k/2}, 
\end{equation} 
where $\xi (\lambda) >1$ and 
\[ \xi (\lambda) = \lambda (1+ O(\lambda^{-1})) \]
in the large coupling regime. 
\end{proposition}
\it Remark. \rm We do not attempt to get the optimal range of values of $\lambda.$ 
Instead, we opt for the clarity of exposition and freely assume that $\lambda$ is large enough. 
The arguments presented below lead to nontrivial dynamical bounds for $\lambda >8.$ With 
more technical effort, this value can be reduced, but remains far from zero. 

\begin{figure}\label{Fig1}
\begin{center}
\setlength{\unitlength}{0.175in}
\begin{picture}(20,7)(0,0)
\put(0,0){\line(1,0){8}} \put(10,0){\hbox to 0mm{\hss$k-1$\hss}}
\put(12,0){\line(1,0){8}}
\put(0,2){\line(1,0){8}} \put(10,2){\hbox to 0mm{\hss$k$\hss}}
\put(12,2){\line(1,0){8}}
\put(0,4){\line(1,0){8}} \put(10,4){\hbox to 0mm{\hss$k+1$\hss}}
\put(12,4){\line(1,0){8}}
\put(0,6){\line(1,0){8}} \put(10,6){\hbox to 0mm{\hss$k+2$\hss}}
\put(12,6){\line(1,0){8}}
\linethickness{2pt}
\put(1,0){\line(1,0){6}}
\put(2,2){\line(1,0){4}}
\put(3,6){\line(1,0){2}}
\put(12.5,2){\line(1,0){7}}
\put(14.5,4){\line(1,0){3}}
\put(13,6){\line(1,0){1}}
\put(18,6){\line(1,0){1}}
\end{picture}
\end{center}
\caption{Types of bands in $\sigma_k:$ left, a type A band; right, a type B band.}
\end{figure}
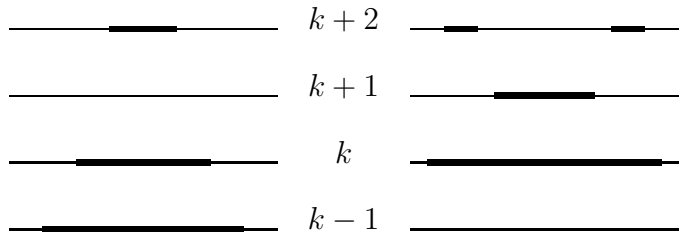

\begin{definition}
We call a band $I_k \subset \sigma_k$ a type A band if $I_k \subset \sigma_{k-1}$ 
(and consequently by Proposition~\ref{background} $I_k \cap (\sigma_{k+1}\cup \sigma_{k-2}) 
= \emptyset$). We call a band $I_k \subset \sigma_k$ a type B band if $I_k \subset \sigma_{k-2}$   
(and so $I_k \cap \sigma_{k-1} = \emptyset$). 
\end{definition}
See Figure~1 for an illustration.
Notice that by the definition of $\Phi_{-1}$ and $\Phi_0,$ 
$\sigma_{-1} = \reals,$ $\sigma_0 =[-2,2]$ and $\sigma_1=[\lambda-2,\lambda+2].$ 
Hence, $\sigma_0$ consists of one band of type A, and $\sigma_1$ consists of one
band of type B. The following Lemma allows to determine inductively the structure 
of the spectrum of $H_\lambda.$

\begin{lemma}\label{struct}
Assume that $\lambda>4.$ Then for any $k \geq 0:$
\begin{enumerate}

\item[i)] Every type A band $I_k \subset \sigma_k$ contains exactly one type B band 
$I_{k+2} \subset \sigma_{k+2}$, and no other bands from $\sigma_{k+1},$ $\sigma_{k+2}.$ 

\item[ii)] Every type B band $I_k \subset \sigma_k$ contains exactly one type A band 
$I_{k+1} \subset \sigma_{k+1}$ and two type B bands from $\sigma_{k+2},$ positioned around 
$I_{k+1}.$ 
\end{enumerate} 
\end{lemma} 
\begin{proof}
Consider a type A band $I_k \subset \sigma_k.$ By definition, $I_k \subset \sigma_{k-1},$ 
and so $|x_{k+1}|>2$ on $I_k$ by Proposition~\ref{background} iii).  This implies 
$I_k \cap \sigma_{k+1} = \emptyset.$
On $I_k,$ $x_k$ changes monotonically from $-2$ to $2,$ in particular, there is a 
unique $E_1 \in I_k$ such that $x_k(E_1)=0.$ By the trace map 
\eqref{TM}, we have $|x_{k+2}(E_1)|=|x_{k-1}(E_1)| \leq 2,$
so $I_k \cap \sigma_{k+2} \ne \emptyset.$ Notice also that when $x_k = \pm 2,$ 
$|x_{k+2}| \geq 2|x_{k+1}|-|x_{k-1}|>2,$ so all possible bands of $\sigma_{k+2}$ intersecting  
$I_k$ lie inside $I_k.$ Moreover, in each band $I_{k+2} \subset I_k$ of $\sigma_{k+2},$ 
$x_{k+2}$ changes from $-2$ 
to $2,$ and so by intermediate value theorem there exists $E_2 \in I_{k+2}$ such that 
$x_{k+2}(E_2) = -x_{k-1}(E_2).$ Then by \eqref{TM}, $x_k(E_2)x_{k+1}(E_2)=0,$ and since 
$|x_{k+1}(E)|>2$ in 
$I_k,$ $x_k(E_2)=0.$ Hence, every band $I_{k+2} \subset I_k$ contains an energy where $x_k=0;$ 
by monotonicity, there is only one such band. 

Now consider a type B band $I_k \subset \sigma_{k-2}.$ When $x_k(E)=0,$ we have
$|x_{k+1}|=|x_{k-2}| \leq 2,$ and so $I_k \cap \sigma_{k+1} \ne \emptyset.$
As in the above argument, $|x_{k+1}|>2$ when $|x_k|= 2,$ so all bands 
of $\sigma_{k+1}$ intersecting $I_k$ lie strictly inside. Moreover, by same argument
as above, any band $I_{k+1} \subset \sigma_{k+1}$ inside $I_k$ must contain an energy 
$E$ where $x_k(E)=0.$ Hence, there is a unique such band $I_{k+1} \subset I_k.$ 
Next consider $\sigma_{k+2}.$ Iterating the trace map, we find
$x_{k+2}=(x_k^2-1)x_{k-1} - x_{k-2}x_k.$ When $x_k = \pm 1,$ $|x_{k+2}|<2.$ Also, if
$x_k =\pm 2,$ $|x_{k+2}| \geq 3|x_{k-1}| -2|x_{k-2}| >2.$ Therefore, there are at least 
two bands of $\sigma_{k+2}$ lying strictly inside $I_k$ to the right and to the left from 
$I_{k+1}$ (bands of $\sigma_{k+2} \subset I_k$ cannot intersect $I_{k+1}$ by iii) of 
Proposition~\ref{background}). 
It remains to show that there are only two such bands. Using \eqref{TM}, it is easy 
to verify that 
\begin{equation}\label{auxTM}
 (x_k \pm 1)(x_{k+2} \pm x_{k-2}) = (x_k^2 -1)(x_{k+1} \pm x_{k-1}). 
\end{equation}
Consider a band $I_{k+2} \subset \sigma_{k+2}$ lying in $I_k.$ For $E \in I_{k+2},$ 
$x_{k+1}(E),$ $x_{k-1}(E)$ 
have fixed signs (both satisfy $|x_{k \pm 1}(E)|>2$ there).    
Pick a sign in \eqref{auxTM} so that $x_{k+1} \pm x_{k-1} \ne 0.$ By the intermediate value 
theorem, there exists the energy 
$E \in I_{k+2}$ where $x_{k+2} \pm x_{k-2} =0.$ At this energy, we must have $x_k^2 -1 =0.$ 
But there are only two energies in $I_k$ where $x_k(E) = \pm 1.$ Hence, there are at most two 
bands of $\sigma_{k+2}$ in $I_k.$   
\end{proof} 

For the proof of the key lemma, we need the following auxiliary result.

\begin{lemma}\label{auxlem}
Let the functions $f_{\pm}(x,y,\lambda)$ be defined by 
\begin{equation}\label{ffun}
f_\pm (x,y, \lambda) = \tfrac{1}{2}\Big[ xy \pm \sqrt{4\lambda^2 + (4-x^2)(4-y^2)}\Big].
\end{equation}
Assume that $\lambda \geq 4.$ 
Then for $-2 \leq x,y \leq 2,$ 
\begin{equation}\label{fbound}
\left| \frac{\prt f_\pm}{\prt x}(x,y,\lambda) \right|,\left| \frac{\prt f_\pm}{\prt y}(x,y,\lambda)
 \right| \leq 1.
\end{equation}
\end{lemma} 
\begin{proof}
Since $f_- (x,y,\lambda)=-f_+(x,-y, \lambda),$ and $f_+(x,y,\lambda)=f_+(y,x,\lambda),$ 
it suffices to show the bound for
\[ \frac{\prt f_+}{\prt x}(x,y,\lambda)= \frac{1}{2}\left(y- \frac{x(4-y^2)}{\sqrt{4\lambda^2
+(4-x^2)(4-y^2)}}\right). \]
Notice that 
\[ \frac{\prt^2 f_+}{\prt x \prt y}(x,y,\lambda) = \frac{1}{2} \left( 1+ 
\frac{8xy \lambda^2 + xy (4-x^2)(4-y^2)}{(4\lambda^2 +(4-x^2)(4-y^2))^{3/2}} \right) \geq 0 \]
for $\lambda >4$ and $|x|,|y| \leq 2.$ 
Thus, 
\[ \max_{|x|,|y| \leq 2} \left| \frac{\prt f_+}{\prt x}(x,y,\lambda) \right| = 
\max_{|x|\leq 2, y=\pm 2} \left| \frac{\prt f_+}{\prt x}(x,y,\lambda) \right| =1. \]
\end{proof}

The Proposition~\ref{tracebound} follows immediately from Lemma~\ref{struct} and

\begin{lemma}\label{keyderlem}
Assume that the coupling $\lambda$ is sufficiently large. Then there exists a number $\xi (\lambda)>1$
such that, given any (type-A) band $I_{k+1} \subset \sigma_{k+1}$ lying in the band 
$I_k \subset \sigma_k,$ we have
\[
\left| \frac{x_{k+1}'(E)}{x_k'(E)} \right| \geq \xi(\lambda)
\]
for $E \in I_{k+1}.$ Similarly, for any (type-B) band $I_{k+2} \subset \sigma_{k+2}$ lying in $I_k,$ 
\[
\left| \frac{x_{k+2}'(E)}{x_k'(E)} \right| \geq \xi(\lambda)
\]
for $E \in I_{k+2}.$ Moreover, $\xi(\lambda) = \lambda (1+O(\lambda^{-1}))$ in the 
large coupling regime. 
\end{lemma}
\begin{proof}
The proof is by induction. The induction is well-founded since $x'_{-1}=0,$ $x'_0=1,$ $x'_1=1.$
There are three cases to consider.  \\
1. First, consider the case where $I_{k+1} \subset \sigma_{k+1}$ lies in $I_k \subset \sigma_k.$ 
Then by Proposition~\ref{background} $I_k \cap \sigma_{k-1}= \emptyset,$ and $I_k \subset \sigma_{k-2}.$ 
Differentiating the equality $x_{k+1} = x_k x_{k-1} - x_{k-2}$ and dividing by $x_k',$ we find
\begin{equation}\label{derrat1}
\frac{x'_{k+1}}{x'_k} = x_{k-1} +\frac{x_k x'_{k-1}}{x'_k}-\frac{x'_{k-2}}{x'_k}. 
\end{equation}
{From} the trace invariant \eqref{Invar}
\[ x_{k}^2+x_{k-1}^2+x_{k-2}^2 -x_{k}x_{k-1} x_{k-2} = 4 +\lambda^2 \]
 we get 
\begin{eqnarray*} 
x_{k-1} = f_\pm(x_k, x_{k-2}, \lambda)=
\tfrac{1}{2} \left( x_k x_{k-2} \pm \sqrt{4\lambda^2 +(4-x_{k-2}^2)(4-x_k^2)}\right) \geq \lambda-2
\end{eqnarray*}
since $|x_k|,|x_{k-2}| \leq 2.$ 
Therefore this term in \eqref{derrat1} is large. To estimate the remaining terms, notice that 
\[ x_{k-1}'= \frac{\prt f_\pm}{\prt x}(x_k, x_{k-2}, \lambda) x'_k + 
\frac{\prt f_\pm}{\prt y}(x_k, x_{k-2}, \lambda) x'_{k-2} \]
(where $\pm$ means that either plus or minus may occur). 
Applying Lemma~\ref{auxlem}, we derive from \eqref{derrat1} that
\[ \left| \frac{x'_{k+1}}{x'_k} \right| \geq \lambda -2 - 3\left| \frac{x'_{k-2}}{x'_k} \right|-|x_k|.\]
By the induction assumption, $|x'_{k-2}/x'_k| \leq \xi(\lambda)^{-1},$ so the induction step holds true 
provided that 
\[ \xi (\lambda) \leq \lambda - 4 - 3\xi(\lambda)^{-1}. \]
The maximal $\xi(\lambda)$ we can take under this condition is 
\[ \xi(\lambda) = \tfrac{1}{2} \left( (\lambda -4) + \sqrt{(\lambda-4)^2 -12} \right) = \lambda +O(1). \]

Next, consider $I_{k+2} \subset \sigma_{k+2}$ lying in $I_k.$ Here we have to distinguish two scenarios. \\
2. $I_{k+2} \cap \sigma_{k-1} = \emptyset,$ $I_{k+2} \subset \sigma_{k-2}.$ 
Given 
\begin{eqnarray*}
x'_{k+2}= x_{k+1}^{} x'_k + x'_{k+1}x_k^{} - x'_{k-1}, \\
x'_{k+1}= x_k^{} x'_{k-1} + x'_k x_{k-1}^{} -x'_{k-2}
\end{eqnarray*}
we find 
\[ \frac{x'_{k+2}}{x'_k} = 2x_{k+1}-x_{k-2}+(x_k^2-1) \frac{x'_{k-1}}{x'_k} -x_k \frac{x'_{k-2}}{x'_k}. \]
Similarly to the previous argument, 
\[ |x_{k+1}| = |f_\pm (x_{k+2}, x_k, \lambda)| \geq \lambda -2. \]
Also, 
\[ x_{k-1}'= \frac{\prt f_\pm}{\prt x}(x_k, x_{k-2}, \lambda) x'_k + 
\frac{\prt f_\pm}{\prt y}(x_k, x_{k-2}, \lambda) x'_{k-2}. \]
Using Lemma~\ref{auxlem}, the bounds $|x_k|,|x_{k-2}| \leq 2,$ and the induction assumption, 
we arrive at 
\[ \left| \frac{x'_{k+2}}{x'_k} \right| \geq 2\lambda -9 - 5\xi(\lambda)^{-1}.\]
Solving the quadratic inequality, we find that the maximal $\xi (\lambda)$ for which induction step
goes through in this case is 
\[ \xi(\lambda) = \tfrac{1}{2} \left( (2\lambda -9) + \sqrt{(2\lambda - 9)^2 -20} \right) = 
2\lambda +O(1). \]
3. $I_{k+2} \subset \sigma_{k-1},$ $I_{k+2} \cap \sigma_{k-2} = \emptyset.$
In this case we obtain 
\[  \frac{x'_{k+2}}{x'_k} = x_{k+1}+ x_k \frac{x'_{k+1}}{x'_k} - \frac{x'_{k-1}}{x'_k}. \]
As before, $x_{k+1} \geq \lambda -2.$ Also 
\[ x'_{k+1} = \frac{\prt f_\pm}{\prt x}(x_k, x_{k-1}, \lambda) x'_k + 
\frac{\prt f_\pm}{\prt y}(x_k, x_{k-1}, \lambda) x'_{k-1}. \]
This leads to 
\[ \left| \frac{x'_{k+2}}{x'_k} \right| \geq \lambda -4 - 3\xi(\lambda)^{-1}.\]
Hence the induction step works with the same $\xi(\lambda)$ as in the first case. 

{From} Lemma~\ref{struct} it follows that we have considered all possible situations.

\end{proof}
\it Remark. \rm In particular, the given argument shows that 
Lemma~\ref{keyderlem} (and hence Proposition~\ref{tracebound}) holds with $\xi(\lambda)>1$ for 
$\lambda \geq 8.$ 

\section{The Fibonacci Hamiltonian: dynamical bounds}

Our first goal is to relate the derivatives of traces with the information required by the
dynamical criterion of Theorem~\ref{lindyn}. 
The following formula is well-known (see, e.g., \cite{Toda})
\begin{align*}
 x'_k(E) =& \sum_{n=1}^{q_k} \Big\{ u_{\pi/2}(q_k,E)u_0^2(n,E) - u_0(q_k+1,E)u_{\pi/2}^2(n,E)+\\
          & \phantom{\sum} + \big[u_{\pi/2}(q_k+1,E)-u_0(q_k,E)\big]u_0(n,E)u_{\pi/2}(n,E) \Big\}
\end{align*}
and can be derived using variation of parameters.
Therefore, 
\[ x'_k(E) \leq 
4 \|\Phi(n,E)\|^3_{q_k+1} \]
(one can show that $x'_k(E) \leq 4(\lambda+O(1))\|\Phi(n,E)\|^2_{q_k+1}$ for $E$ in the spectrum 
of $H_\lambda,$ but we choose to avoid the technicalities).
Hence, Proposition~\ref{tracebound} implies 
\begin{equation}\label{eqf1}
 \|\Phi(n,E)\|_{q_k+1}^3 \geq \frac{1}{4} \xi(\lambda)^{k/2}. 
\end{equation}
As $\omega$ is the golden ratio $(\sqrt{5}-1)/2$,  
$$
q_k = [\omega^{-k} - (-\omega)^{k}]/\sqrt{5}
$$
for $k \geq 1$ (see, e.g., \cite{IT1}).
So for large $k$, $q_k \sim \omega^{-k}/\sqrt{5}$. Consider $q_k \leq L \leq q_{k+1}$.
Since $q_{k+1}/q_k \leq 2$ for all $k$, \eqref{eqf1} implies
\begin{equation}\label{tranbound}
 \| \Phi(n,E)\|^2_L \geq C q_k^{\zeta_1} \geq C' L^{\zeta_1} \quad \text{ with } 
\zeta_1= {\frac{\log \xi(\lambda)}{3 \log (\omega^{-1})}}.
\end{equation}
Later we will derive an upper bound on this quantity of the form $CL^{\zeta_2}$ with $\zeta_2$ 
also depending logarithmically on $\lambda$.

Now we complete the proof of the upper dynamical bound in Theorem~\ref{mainfib}. 

\begin{proof}
By the symmetry of Fibonacci potential ($V(-n)=V(n-1)$ for $n \geq 2$, see \cite{Suto1}) the bound
identical to \eqref{tranbound} also holds for the negative semi-axis. These bounds allow us to define 
the characteristic scales $\tilde{L}^\pm_{T^{-1}}(E).$  
Consider the case of $\tilde{L}^+_{T^{-1}}(E)$ (the other one is analogous). 
By \eqref{chartran}, $\tilde{L}^+_{T^{-1}}(E)$ is determined by the equality
\[ \|\Phi(n, E)\|^2_{\tilde{L}^+_{T^{-1}}(E)} =4 \|\Phi(1,E)^{-1}\|^2  T^2 .\]
{From} \eqref{tranbound} we get 
\[ \tilde{L}^+_{T^{-1}}(E)^{\frac{\log \xi(\lambda)}{3\log (\omega^{-1})}} \leq C T^2, \]
and since $\xi(\lambda) = \lambda + O(1)$ we can take 
\[ \tilde{L}^+_{T^{-1}}(E) = C T^\frac{6 \log (\omega^{-1})}{\log \xi(\lambda)}. \]
This proves the first part of Theorem~\ref{mainfib}. 
\end{proof}

It remains to prove the second part of Theorem~\ref{mainfib}, involving the lower bound on dynamics. 
The lower bound on dynamics for the Fibonacci Hamiltonian has been proved recently in \cite{JL2}.
We present a sketch of the argument here, making explicit the behavior of the bound in the 
large coupling regime. The lower dynamical bound will follow from the continuity estimate for
the spectrum. 
The idea is to study the behavior of solutions; as soon as appropriate bounds are available, 
one could apply the reasoning of Jitomirskaya-Last \cite{JL2},  
a recent result of Damanik, Killip and Lenz \cite{DKL}
relating bounds on solutions to the $\alpha$-continuity of the whole line operator,
or Theorem~\ref{ldb}.

\begin{lemma}\label{unifbeh}
For every energy $E$ in the spectrum of $H_\lambda,$ every solution $u(n,E)$ of the equation 
$(H_\lambda-E)u(n,E)=0$ satisfies 
\begin{equation}\label{unifsol}
\|u(n,E)\|_L \geq CL^{\kappa}, 
\end{equation}
where $\kappa = \log \frac{\sqrt{17}}{20 \log (\omega^{-1})}$ is independent of $\lambda.$ 
\end{lemma} 
\begin{proof}
The proof follows closely Proposition 10 of \cite{JL2}. 
As shown by S\"ut\H{o} \cite{Suto1}, the Fibonacci potential obeys $V(q_n+l)=V(l)$ for $n \geq 3$ 
and $1 \leq l \leq q_n.$ For $n$ large and $1 \leq l \leq q_{n-2},$ we also get 
$V(2q_n+l)= V(q_{n+1}+q_{n-2}+l)=V(q_{n-2}+l)=V(l),$ so that 
\begin{equation}\label{square}
V(l)= V(q_n +l)= V(2q_n +l).
\end{equation} 
By Lemma 1 of \cite{Suto1}, we have for any $2 \times 2$ matrix $B$ with ${ \rm det}B=1,$ 
\[ {\rm max}\{ |\tr B| \|B \Psi\|, \|B^2 \Psi \| \} \geq \frac{1}{2}\|\Psi\| \]
 for any $2$-vector $\Psi.$ Therefore, 
\begin{equation}\label{corsq}
\|B\Psi\|^2 + \|B^2 \Psi\|^2 > \frac{1}{4{\rm max}(1,|\tr B|^2)} \|\Psi\|^2. 
\end{equation}
Let $\Phi(m,k)$ denote the transfer matrix that takes 
$(u(m+1),u(m))^T$ to $(u(k+1),u(k))^T.$ By \eqref{square}, for $1 \leq l \leq q_{n-2}$ we have 
$\Phi(l, q_n+l)=\Phi(q_n+l,2q_n+l),$ and, moreover, $\tr \Phi (l, q_n+l)=x_n.$  
Thus, \eqref{corsq} implies
\begin{eqnarray*}
 |u(q_n+l+1)|^2 + |u(q_n +l)|^2 + |u(2q_n +l+1)|^2 +|u(2q_n+l)|^2>  \\ \frac{1}{4{\rm max}(1,x_n^2)}
(|u(l+1)|^2+|u(l)|^2) 
\end{eqnarray*}
for $1 \leq l \leq q_{n-2},$ and 
\begin{eqnarray*}
 |u(q_{n+1}+l+1)|^2 + |u(q_{n+1} +l)|^2 + |u(2q_{n+1} +l+1)|^2 +|u(2q_{n+1}+l)|^2> \\ 
\frac{1}{4{\rm max}(1,x_{n+1}^2)}
(|u(l+1)|^2+|u(l)|^2) 
\end{eqnarray*}
for $1 \leq l \leq q_{n-1}.$
Notice that for any $E$ in the spectrum of $H_\lambda,$ for every $n$ either $|x_n|$ or 
$|x_{n+1}|$ is less than $2.$ 
By combining the above estimates, we easily deduce that 
\begin{equation}\label{keyit}
 \|u(n,E)\|_{q_{n+5}} \geq \frac{\sqrt{17}}{4} \|u(n,E)\|_{q_n} 
\end{equation}
for large $n.$ The estimate \eqref{unifsol} follows directly from \eqref{keyit}.  
\end{proof}

On the other hand, Iochum and Testard \cite{IT1} have shown that 
\[ \|\Phi(n,E)\| \leq \left( d^{\frac{\log \sqrt{5}}{\log (\omega^{-1})}} \right)^{\log n}. \]
The number $d$ is defined as a product $d=ab^2,$ where $b=2c+1,$ $c={\rm sup}_n |x_n|,$ and
$a={\rm max} (c,2).$ For $\lambda$ large, one can check using \eqref{Invar} that 
$c \leq \lambda +2$.
Hence, $d \leq (\lambda+2)(2\lambda +5)^2.$ Therefore we conclude that 
\begin{equation}\label{tranupbound}
\sum\limits_{n=1}^L \|\Phi(n,E)\|^2 \leq CL^{\zeta_2(\lambda)}, 
\end{equation}
where 
\[ \zeta_2(\lambda)= \frac{6 \log \sqrt{5}}{\log (\omega^{-1})}(\log \lambda +O(1)). \]

Recall that a measure $\mu$ is called $\alpha$-continuous if it does not give weight to sets 
of zero $\alpha$-dimensional Hausdorff measure.

\begin{proposition}\label{acont}
The spectral measure $\mu$ of the operator $H_\lambda$ is $\alpha(\lambda)$-continuous
for every $\lambda>0,$ with $\alpha(\lambda)$ satisfying 
\begin{equation}\label{alpha}
\alpha(\lambda) \geq \frac{2\kappa}{\kappa +\zeta_2(\lambda)}=C (\log \lambda)^{-1}(1+O(( \log \lambda)^{-1})). 
\end{equation}
\end{proposition}
\begin{proof}
Given the estimates \eqref{unifsol} and \eqref{tranupbound}, for the Fibonacci operator $H_\lambda$
the Proposition follows along the lines of \cite{JL2}. Alternatively, one may use 
 Theorem 1 of \cite{DKL} which shows for any discrete 
Schr\"odinger operator $H$ on whole axis
that if every solution of $(H-E)u=0$ obeys
\[ C_1(E)L^{q_1} \leq \|u(n,E)\|_L \leq C_2(E)L^{q_2} \]
for $E \in \Sigma,$
then the spectrum of $H$ is $\alpha$-continuous in $\Sigma$ with $\alpha=2q_1/(q_1+q_2).$ 
\end{proof}   

Now we complete the proof of Theorem~\ref{mainfib}. 
\begin{proof}
The fact that $\alpha$-continuity of the spectral measure implies lower bounds on dynamics is well-known. 
The original result is due to Guarneri \cite{Guarneri,Gu2}, with later contributions by many authors.
Assume that the spectral measure $\mu^\phi$ of the Schr\"odinger operator $H$
corresponding to the vector $\phi$ 
is $\alpha$-continuous. Then by Theorem 4.2 of \cite{Last}, the measure $\mu^\phi$ can be represented
as a sum of mutually singular measures $\mu^\phi = \mu^{\phi_1} +\mu^{\phi_2},$ 
where $\mu^{\phi_1}$ is uniformly $\alpha$-H\"older continuous ($\mu^{\phi_1}(I) \leq C|I|^\alpha$ for any 
interval $I$) while $\mu^{\phi_2}(\reals)$ can be made arbitrarily small.
Moreover, it has been shown in \cite{Last}, proof of Theorem 6.1, that 
\begin{equation}\label{dynboun}
\tfrac{1}{a} \int\limits_0^a \|e^{-iHt}\phi_1\|_L^2 \, dt \leq CL a^{-\alpha}.
\end{equation}
(in dimension one).
Multiplying  \eqref{dynboun} by $a$ and integrating with respect to the kernel $\frac{2}{T^2} e^{-2a/T},$ 
we get 
\[  \tfrac{2}{T} \int_0^\infty e^{-2t/T} \|e^{-iHt}\phi_1\|_L^2 \, dt \leq C' L T^{-\alpha}. \]
In particular, given any $g>0,$ we can choose $\phi_2$ of sufficiently small norm and
$C_g$ so that for $L(T) = C_g T^\alpha,$
\begin{equation}\label{dynfin}
  \tfrac{2}{T} \int_0^\infty e^{-2t/T} \|e^{-iHt}\phi\|_{L(T)}^2 \, dt \leq g. 
\end{equation}
Application of \eqref{dynfin} and Proposition~\ref{acont} completes the proof.
\end{proof}

\section{Multidimensional problems}

Our purpose in this section is to show how our results for tridiagonal operators
can be applied to more general problems. We are particularly interested in discrete
Schr\"odinger operators of the form $\Delta+V$ on $\ell^2(\Z^d)$,
defined by
\begin{equation}\label{mdj}
((\Delta+V)\psi)(n) = \sum_{|n-m|=1}\psi(m)+V(n)\psi(n),
\end{equation}
where $\{V(n)\}_{n\in\Z^d}\subset\R$.

Given any self adjoint operator $H$ on a separable Hilbert space $\hils$ and a vector
$\psi$ in the domain of $H$,
the cyclic subspace spanned by $H$ and $\psi$ is defined by
$\hils_\psi = \overline{\left\{f(H)\psi\left|\right. f\in C_\infty(\R)\right\}}$.
Here $C_\infty(\R)$ is the set of continuous (complex valued) functions on $\R$ vanishing at infinity
and $\overline{\;\,\cdot\;\,}$ denotes norm closure in $\hils$.  
$\hils_\psi$ is an
invariant subspace for $H$ and by the spectral theorem, the restricted operator
$H\restriction\hils_\psi$ is unitarily equivalent to multiplication by the
coordinate parameter on the space $L^2(\R,d\mu_\psi)$, where $\mu_\psi$ is
the spectral measure of $\psi$ (and $H$). We assume
that $\mu_\psi$ is not supported on a finite number of points and thus $\hils_\psi$
is infinite dimensional.
The evolution given by \eqref{evolution} is thus confined to $\hils_\psi$
(namely, $\psi(t)\in\hils_\psi$ for any $t$) and so to study the time evolution,
it suffices to consider the restricted
operator $H\restriction\hils_\psi$.

For large classes of problems, the moment vectors $\{H^n\psi\}_{n=0}^\infty$ are
well defined and the set of their finite linear combinations is dense in $\hils_\psi$,
namely, 
\begin{equation}\label{span}
\hils_\psi = \overline{\text{span}\{H^n\psi\}_{n=0}^\infty} .
\end{equation}
We note that this holds for any $\psi$ if $H$ is bounded. If $H$ is unbounded, then $\psi$
must be chosen appropriately in order for this to occur (but such $\psi$'s
always exist, e.g., any $\psi$ with a compactly supported spectral measure). We note
that $\psi$ must be in the domain of each of the moment operators $H^n$ and thus
its spectral measure must have sufficient decay at infinity to ensure finite moments.
This necessary condition is not sufficient, though, and determining
whether this property holds or not for a given $\psi$ is in general a rather rich question,
equivalent to the problem of whether the spectral measure of $\psi$ can be determined
from its moments. See \cite{Simon3} for more information. The important thing for us here
is that the property \eqref{span} holds whenever $H$ is of the form \eqref{mdj} on
$\ell^2(\Z^d)$ and $\psi$ is any vector which is compactly supported on $\Z^d$.

If the property \eqref{span} holds, then $\hils_\psi$ has a natural orthonormal
basis $\{\rho_n\}_{n=0}^\infty$, which is obtained by applying the Gram-Schmidt
orthonormalization procedure to the set of moment vectors $\{H^n\psi\}_{n=0}^\infty$.
Each $\rho_n$ has the form $\rho_n = P_n(H)\psi$, where $P_n(\cdot)$ is a polynomial of
degree $n$. The polynomials $\{P_n(\cdot)\}_{n=0}^\infty$ are known as the orthonormal
polynomials of the spectral measure $\mu_\psi$. The self-adjoint operator
$H\restriction\hils_\psi$ is known \cite{Simon3} to have a tridiagonal matrix representation in the
basis $\{\rho_n\}_{n=0}^\infty$ and by choosing the phases of the $\rho_n$ appropriately,
this tridiagonal matrix can be made to have only real entries. Moreover, the off-diagonal
entries cannot vanish. Thus, we see that $\hils_\psi$ can be viewed as being
$\ell^2(\Z_+)$, by identifying $\rho_n$ with $\delta_{n+1}$, and then
$H\restriction\hils_\psi$ has precisely the form \eqref{jacobi} with a Dirichlet
boundary condition. Theorems~\ref{hld}, \ref{hldt}, \ref{ldb}, and \ref{pb} are thus
fully applicable in this general setting. We call the above tridiagonal
representation of $H\restriction\hils_\psi$ the
orthonormal polynomial representation.

The problem, of course, is that when one is interested in wavepacket spreading
dynamics, one is normally interested in how the spreading occurs in some natural
coordinate space which is generally different from the orthonormal basis
$\{\rho_n\}_{n=0}^\infty$. The ability to apply our results to such a coordinate
space depends on the ability to relate the vectors $\rho_n$ to it. Fortunately,
for the case of an operator of the form \eqref{mdj} on $\ell^2(\Z^d)$ and a
compactly supported initial vector $\psi$, an appropriate connection exists. For
simplicity, take the initial vector $\psi$ to be the delta function vector at the
origin of $\Z^d$, namely $\psi=\delta_0$. Then one easily sees from \eqref{mdj}
that $H^n\delta_0$, and thus also $\rho_n$, is supported in a ball of radius
$n+1$ in $\Z^d$. Thus, any portion of the wavepacket which is supported
on the vectors $\rho_0,...,\rho_L$ must also be confined to a ball of radius
$L+1$ in $\Z^d$. This means that any upper bound on the spreading of the wavepacket
in the ``coordinates'' $\{\rho_n\}_{n=0}^\infty$ is immediately also bounding
its spreading in $\Z^d$. In particular, Theorem~\ref{hld} takes the following
form.
\begin{theorem}\label{mdhld}
Let $H$ be an operator of the form \eqref{mdj} on $\ell^2(\Z^d)$.
Let the characteristic scale $L^+_{T^{-1}}(E)$ be defined by \eqref{ls}
for the orthonormal polynomial representation of $H\restriction\hils_{\delta_0}$.
Then for any $T>0$ and $L>1$, we have 
\begin{equation}\label{mddynfir}
\langle \|\chi_{L+1} e^{-iH t} \delta_0 \|^2 \rangle_T >
C \mu_{\delta_0}\left(\left\{ E \left| \right. L^+_{T^{-1}}(E) \leq L \right\}\right),
\end{equation} 
where $C$ is some universal positive constant and $\chi_L$ denotes the orthogonal
projection on a ball of radius $L$ around the origin in $\Z^d$.
\end{theorem}
We note that the scale $L^+_{T^{-1}}(E)$ in Theorem \ref{mdhld} is defined through
the orthonormal polynomial representation, a fact which may make it seem like
a somewhat obscured quantity in the multidimensional context. It is interesting
to note, however, that even in the multidimensional context, it is actually among
the more straight forward quantities to compute. Each of the moment vectors
$H^n\delta_0$ can be easily computed explicitly and it only depends on values
of the potential within a ball of a corresponding radius $n$, and given the
moment vectors up to some scale $L$, the corresponding $\rho_0,...,\rho_L$ are
obtained from them by the Gram-Schmidt procedure. The matrix elements of the
(tridiagonal) orthonormal polynomial representation of $H\restriction\hils_{\delta_0}$
are then easily obtained up to a corresponding scale $L$ and from them one
can explicitly compute solutions of equation \eqref{ef} up to a corresponding
scale. It is thus possible, for any $E$ and $T$, to compute (or more precisely,
to bound within two integer values) the scale $L^+_{T^{-1}}(E)$ by a straight
forward computation consisting of finitely many steps and only depending on
values of the potential within a ball of a corresponding radius. This should
be contrasted with the fact that multidimensional representations of generalized
eigenfunctions, such as the ones considered in the Kiselev-Last lower bound \cite{KL},
depend on the entire potential and cannot be computed by a finite
procedure.

\smallskip\bigskip
\noindent {\bf Acknowledgment.} 
We would like to thank Svetlana Jitomirskaya and Barry Simon for useful discussions.

\end{document}